
\def\ignore#1{}
 

\newcount\sectnum
\newcount\subsectnum
\newcount\eqnumber

\global\eqnumber=1\sectnum=0


\def\lab{(\the\sectnum.\the\eqnumber)}



\def\show#1{#1}



\def\smskip{\vskip 5 pt}
\def\medskip{\vskip 10 pt}
\def\bigskip{\vskip 15 pt}
\def\pn{\par\noindent}
\def\br{\break}

\def\bl{\bigl} 
\def\br{\bigr} 
\def\lf{\left}
\def\ri{\right}

\def\argmin{\mathop{\arg \min}}

\def\frac#1#2{{#1\over #2}}

\def\ubar{\underline}
\def\ol#1{\overline{#1}}

\def\a{\alpha}

\def\b{\beta}
\def\l{\lambda}
\def\g{\gamma}
\def\m{\mu}
\def\p{\pi}
\def\r{\rho}
\def\e{\epsilon}

\def\d{\delta}

\def\P{\Pi}

\def\re{\Re}
\def\rn{\Re^n}

\def\tl{\tilde}

\def\old#1{}
\def\leaderfill{\leaders\hbox to 1em{\hss.\hss}\hfill}


\parindent=2pc
\baselineskip=15pt
\vsize=8.7 true in
\voffset=0.125 true in
\parskip=3pt


\def\minprob#1#2#3{$$\eqalign{&\hbox{minimize\ \ }#1\cr &\hbox{subject to\ \
}#2\cr}\ifnum 0=#3{}\else\eqno(#3)\fi$$}        
     
\def\maxprob#1#2#3{$$\eqalign{&\hbox{maximize\ \ }#1\cr &\hbox{subject to\ \
}#2\cr}\ifnum 0=#3{}\else\eqno(#3)\fi$$}        
     
\def\aligntwo#1#2#3#4#5{$$\eqalign{#1&#2\cr #3&#4\cr}
\ifnum 0=#5{}\else\eqno(#5)\fi$$}
\def\alignthree#1#2#3#4#5#6#7{$$\eqalign{#1&#2\cr #3&#4\cr #5&#6\cr}
\ifnum 0=#7{}\else\eqno(#7)\fi$$}


\def\eqnum{\eqno{\hbox{(\the\sectnum.\the\eqnumber)}\global\advance\eqnumber
by1}}

\def\eqnu{\eqno{\hbox{(\the\sectnum.\the\eqnumber)}\global\advance\eqnumber
by1}}

\newcount\examplnumber
\def\examplnum{\global\advance\examplnumber by1}

\newcount\figrnumber
\def\figrnum{\global\advance\figrnumber by1}

\newcount\propnumber
\def\propnum{\global\advance\propnumber by1}

\newcount\defnumber
\def\defnum{\global\advance\defnumber by1}

\newcount\lemmanumber
\def\lemmanum{\global\advance\lemmanumber by1}

\newcount\assumptionnumber
\def\assumptionnum{\global\advance\assumptionnumber by1}

\newcount\conditionnumber
\def\conditionnum{\global\advance\conditionnumber by1}

\def\exampl{\the\sectnum.\the\examplnumber}
\def\figr{\the\sectnum.\the\figrnumber}
\def\propn{\the\sectnum.\the\propnumber}
\def\defn{\the\sectnum.\the\defnumber}
\def\lemman{\the\sectnum.\the\lemmanumber}
\def\assumptionn{\the\sectnum.\the\assumptionnumber}
\def\condn{\the\sectnum.\the\conditionnumber}

\def\section#1{\goodbreak\vskip 3pc plus 6pt minus 3pt\leftskip=-2pc
   \global\advance\sectnum by 1\eqnumber=1
\global\examplnumber=1\figrnumber=1\propnumber=1\defnumber=1\lemmanumber=1\assumptionnumber=1 \conditionnumber =1%
   \line{\hfuzz=1pc{\hbox to 3pc{\bf 
   \vtop{\hfuzz=1pc\hsize=38pc\hyphenpenalty=10000\noindent\uppercase{\the\sectnum.\quad #1}}\hss}}
			\hfill}
			\leftskip=0pc\nobreak\tenf
			\vskip 1pc plus 4pt minus 2pt\noindent\ignorespaces}



\def\sect#1{\noindent\leftskip=-2pc\tenf
   \goodbreak\vskip 1pc plus 4pt minus 2pt
                \global\advance\subsectnum by 1\eqnumber=1
   \line{\hfuzz=1pc{\hbox to 3pc{\bf 
   \vtop{\hfuzz=1pc\hsize=38pc\hyphenpenalty=10000\noindent\uppercase{{\bf #1}}}\hss}}
                        \hfill}
   \leftskip=0pc\nobreak\tenf
                        \vskip 1pc plus 4pt minus 2pt\nobreak\noindent\ignorespaces}

\def\subsection#1{\noindent\leftskip=0pc\tenf
   \goodbreak\vskip 1pc plus 4pt minus 2pt
   \line{\hfuzz=1pc{\hbox to 3pc{\bf 
   \vtop{\hfuzz=1pc\hsize=38pc\hyphenpenalty=10000\noindent{\bf #1}}\hss}}
                        \hfill}
   \leftskip=0pc\nobreak\tenf
                        \vskip 1pc plus 4pt minus 2pt\nobreak\noindent\ignorespaces}
\def\subsubsection#1{\goodbreak\vskip 1pc plus 4pt minus 2pt
   \hfuzz=3pc\leftskip=0pc\noindent\tenit #1 \nobreak\tenf\vskip 6pt plus 1pt
                                minus 1pt\nobreak\ignorespaces\leftskip=0pc}
%

\def\beginexample#1{\noindent\goodbreak\vskip 6pt plus 1pt minus 1pt
\noindent
  \hbox {\bf Example #1\hss}
  \nobreak\vskip 4pt plus 1pt minus 1pt \nobreak\noindent\ninef
  \global\advance
                \leftskip by\parindent\pn}
\def\endexample{\vskip 12pt\tenf\par
  \global\advance\leftskip by -\parindent
  }

\def\beginexercise#1{\noindent\goodbreak\vskip 6pt plus 1pt minus 1pt \noindent\global\normalbaselineskip=12pt
  \hbox {\bf Exercise #1\hss}
  \nobreak\vskip 4pt plus 1pt minus 1pt 
  \nobreak\noindent\ninef\global\advance\leftskip
                        by\parindent\pn}
\def\endexercise{\vskip 12pt\tenf\par
  \global\advance\leftskip by -\parindent
  }

\def\beginsection#1{\noindent\goodbreak\vskip 6pt plus 1pt minus 1pt \noindent\global\normalbaselineskip=12pt
  \hbox {\it #1\hss}
  \vskip 0.1pt plus 1pt minus 1pt \nobreak\noindent\ninef\global\advance
                \leftskip by\parindent\noindent\pn}
\def\endsection{\vskip 12pt\tenf\par
  \global\advance\leftskip by -\parindent
}

%


\def\proposition#1{\smskip\pn{\bf Proposition #1}\quad}
\def\proof{\smskip\pn{\bf Proof:}\quad} 
\def\definition#1{\smskip\pn{\bf Definition #1}\quad} 
\def\assumption#1{\smskip\pn{\bf Assumption #1}\quad}

 \def\qed{\quad{\bf
Q.E.D.} \par\bigskip}
\def\ref{\smskip\pn}

\def\chapter#1#2{{\bf \centerline{\helbigbig
{#1}}}\bigskip\bigskip{\bf \centerline{\helbigbig
{#2}}}\bigskip\bigskip} 



\def\longpapertitle#1#2#3{{\bf \centerline{\helbigb
{#1}}}\bigskip{\bf \centerline{\helbigb
{#2}}}\bigskip\bigskip{\centerline{
by}}\bigskip{\bf \centerline{
{#3}}}\bigskip\bigskip} 


\def\nitem#1{\smskip\item{#1}}

\newcount\alphanum
\newcount\romnum

\def\alphaenumerate{\ifcase\alphanum \or (a)\or (b)\or (c)\or (d)\or (e)\or
(f)\or (g)\or (h)\or (i)\or (j)\or (k)\fi}
\def\romenumerate{\ifcase\romnum \or (i)\or (ii)\or (iii)\or (iv)\or (v)\or
(vi)\or (vii)\or (viii)\or (ix)\or (x)\or (xi)\fi}

\def\alist{\begingroup\vskip10pt\alphanum=1
\parskip=2pt\parindent=0pt \leftskip=3pc
\everypar{\llap{{\rm\alphaenumerate\hskip1em}}\advance\alphanum by1}}

\def\nolist{\begingroup\vskip10pt\alphanum=0
\parskip=2pt\parindent=0pt \leftskip=3pc
\everypar{\llap{\global\advance\alphanum by1(\the\alphanum)\hskip1em}}}

\def\romlist{\begingroup\vskip10pt\romnum=1
\parskip=2pt\parindent=0pt \leftskip=5pc
\everypar{\llap{{\rm\romenumerate\hskip1em}}\advance\romnum by1}}



\long\def\fig#1#2#3{\vbox{\vskip1pc\vskip#1
\prevdepth=12pt \baselineskip=12pt
\vskip1pc
\hbox to\hsize{\hfill\vtop{\hsize=25pc\noindent{\eightbf Figure #2\ }
{\eightpoint#3}}\hfill}}}

\long\def\widefig#1#2#3{\vbox{\vskip1pc\vskip#1
\prevdepth=12pt \baselineskip=12pt
\vskip1pc
\hbox to\hsize{\hfill\vtop{\hsize=28pc\noindent{\eightbf Figure #2\ }
{\eightpoint#3}}\hfill}}}

\long\def\table#1#2{\vbox{\vskip0.5pc
\prevdepth=12pt \baselineskip=12pt
\hbox to\hsize{\hfill\vtop{\hsize=25pc\noindent{\eightbf Table #1\ }
{\eightpoint#2}}\hfill}}}

 
\def\rightheadline#1{\headline{\tenrm\hfil #1}}


\long\def\leftfig#1#2{\vbox{\smskip\hsize=220pt
\vtop{{\noindent {\bf #1}}}
\smskip
\noindent
\vbox{{\noindent #2}}
}}

\long\def\rightfig#1#2#3{\vbox{\smskip\vskip#1
\prevdepth=12pt \baselineskip=12pt
\hsize=210pt
\smskip
\vbox{\noindent{\eightbold #2}
\hskip1em{\eightpoint#3}}
}}

\long\def\concept#1#2#3#4#5{\bigskip\hrule
\vbox{\hbox{\leftfig{#1}{#2} \hskip3em
\rightfig{#3}{#4}{#5}} \smskip}
\hrule\bigskip}


\long\def\bconcept#1#2#3#4#5#6#7{
\vbox{
\hbox to \hsize{\vtop{\par #1}}
\concept{#2}{#3}{#4}{#5}{#6}
\hbox to \hsize{\vtop{\par #7}}
\smskip}
}




\def\boxit#1{\vbox{\hrule\hbox{\vrule\kern3pt
                                \vbox{\kern3pt#1\kern3pt}\kern3pt\vrule}\hrule}}
\def\centerboxit#1{$$\vbox{\hrule\hbox{\vrule\kern3pt
                                \vbox{\kern3pt#1\kern3pt}\kern3pt\vrule}\hrule}$$}

\long\def\boxtext#1#2{$$\boxit{\vbox{\hsize #1\noindent\strut #2\strut}}$$}

%
%
%

\def\picture #1 by #2 (#3){
  \vbox to #2{
    \hrule width #1 height 0pt depth 0pt
    \vfill
    \special{picture #3} 
    }
  }

\def\scaledpicture #1 by #2 (#3 scaled #4){{
  \dimen0=#1 \dimen1=#2
  \divide\dimen0 by 1000 \multiply\dimen0 by #4
  \divide\dimen1 by 1000 \multiply\dimen1 by #4
  \picture \dimen0 by \dimen1 (#3 scaled #4)}
  }

%
%

\long\def\captfig#1#2#3#4#5{\vbox{\vskip1pc
\hbox to\hsize{\hfill{\picture #1 by #2 (#3)}\hfill}
\prevdepth=9pt \baselineskip=9pt
\vskip1pc
\hbox to\hsize{\hfill\vtop{\hsize=24pc\noindent{\eightbold Figure #4}
\hskip1em{\eightpoint#5}}\hfill}}}

%
%
%

\def\illustration #1 by #2 (#3){
  \vskip#2\hskip#1\special{illustration #3} 
    }

\def\scaledillustration #1 by #2 (#3 scaled #4){{
  \dimen0=#1 \dimen1=#2
  \divide\dimen0 by 1000 \multiply\dimen0 by #4
  \divide\dimen1 by 1000 \multiply\dimen1 by #4
  \illustration \dimen0 by \dimen1 (#3 scaled #4)}
  }


\newbox\graybox
\newdimen\xgrayspace
\newdimen\ygrayspace
%
%
%
%
%
%
%
%
%

\def\Textshade#1#2#3#4#5#6{%
    \xgrayspace=#4pt%
    \ygrayspace=#4pt%
    \def\grayshade{#3}%
    \def\linewidth{#5}%
    \def\theradius{#6}%
    \setbox\graybox=\hbox{\surroundboxa{#2}}%
    \hbox{%
    \hbox to 0pt{%
    \PScommands
}%
    \box\graybox}}%
%
%
\long%

\long%
\def\Parashade#1#2#3#4#5#6#7{%
    \xgrayspace=#4pt%
    \ygrayspace=#4pt%
    \def\grayshade{#3}%
    \def\linewidth{#5}%
    \def\theradius{#6}%
    \def\thevskip{#7pt}%
    \setbox\graybox=\hbox{\surroundboxb{#2}}%
    \vskip\thevskip%
    \hbox{%
    \hbox to 0pt{%
    \PScommands
}%
     \box\graybox}%
     \vskip\thevskip%
}%
%
%
%
\long\def\surroundboxa#1{\leavevmode\hbox{\vtop{%
\vbox{\kern\ygrayspace%
\hbox{\kern\xgrayspace#1%
      \kern\xgrayspace}}\kern\ygrayspace}}}
%
%
\long\def\surroundboxb#1{\leavevmode\hbox{\vtop{%
\vbox{\kern\ygrayspace%
\hbox{\kern\xgrayspace\vbox{\advance\hsize-2\xgrayspace#1}%
      \kern\xgrayspace}}\kern\ygrayspace}}}
%
%
%
\long\def\PScommands{%
\special{rawpostscript
/sharpbox{%
           newpath
           xmin ymin moveto
           xmin ymax lineto
           xmax ymax lineto
           xmax ymin lineto
           xmin ymin lineto
           closepath 
          }bind def
}%
\special{rawpostscript
/sharpboxnb{%
           newpath
           xmin ymin moveto
           xmin ymax lineto
           xmax ymax lineto
           xmax ymin lineto
          }bind def
}%
\special{rawpostscript
/sharpboxnt{%
           newpath
           xmin ymax moveto
           xmin ymin lineto
           xmax ymin lineto
           xmax ymax lineto
          }bind def
}%
\special{rawpostscript
/roundbox{%
           newpath
           xmin radius add ymin moveto
           xmax ymin xmax ymax radius arcto
           xmax ymax xmin ymax radius arcto
           xmin ymax xmin ymin radius arcto
           xmin ymin xmax ymin radius arcto 16 {pop} repeat
           closepath
          }bind def
}%
\special{rawpostscript
/sharpcorners{%
               sharpbox gsave grayshade setgray fill grestore 
               linewidth setlinewidth stroke
              }bind def
}%
\special{rawpostscript
/sharpcornersnt{%
               sharpboxnt gsave grayshade setgray fill grestore 
               linewidth setlinewidth stroke
              }bind def
}%
\special{rawpostscript
/sharpcornersnb{%
               sharpboxnb gsave grayshade setgray fill grestore 
               linewidth setlinewidth stroke
              }bind def
}%
\special{rawpostscript
/roundcorners{%
               roundbox gsave grayshade setgray fill grestore 
               linewidth setlinewidth stroke
              }bind def
}%
\special{rawpostscript
/plainbox{%
           sharpbox grayshade setgray fill 
          }bind def
}%
%
\special{rawpostscript
/roundnoframe{%
               roundbox grayshade setgray fill 
              }bind def
}%
\special{rawpostscript
/sharpnoframe{%
               sharpbox grayshade setgray fill 
              }bind def
}%
}%
%
%

\def\pshade#1{\Parashade{sharpcorners}{#1}{0.95}{10}{0.5}{10}{10}}


\def\boxit#1{\vbox{\hrule\hbox{\vrule\kern3pt
                                \vbox{\kern3pt#1\kern3pt}\kern3pt\vrule}\hrule}}

\def\boxitnb#1{\vbox{\hrule\hbox{\vrule\kern3pt
                                \vbox{\kern3pt#1\kern3pt}\kern3pt\vrule}}}

\def\boxitnt#1{\vbox{\hbox{\vrule\kern3pt
                                \vbox{\kern3pt#1\kern3pt}\kern3pt\vrule}\hrule}}

\long\def\boxtext#1#2{$$\boxit{\vbox{\hsize #1\noindent\strut #2\strut}}$$}



\def\texshopbox#1{\boxtext{462pt}{\vskip-1.5pc\pshade{\vskip-1.0pc#1\vskip-2.0pc}}}


%
%
%
%
%
%
%
%
\font\helbigbig=cmr10 scaled 2500%
\font\helbigb=cmbx10 scaled 1500%
\font\eightbold=cmbx8%

\def\tenf{\hel}%
\def\tenit{\heli}%
\def\ninef{\ninehel}%
\def\nineit{\nineheli}%
%
%


\font\tenrm=cmr10%
\font\teni=cmmi10%
\font\tensy=cmsy10%
\font\tenbf=cmbx10%
\font\tentt=cmtt10%
\font\tenit=cmti10%
\font\tensl=cmsl10%

\def\tenpoint{\def\rm{\fam0\tenrm}%
\textfont0=\tenrm%
\textfont1=\teni%
\textfont2=\tensy%
\textfont\itfam=\tenit%
\textfont\slfam=\tensl%
\textfont\ttfam=\tentt%
\textfont\bffam=\tenbf%
\scriptfont0=\sevenrm%
\scriptfont1=\seveni%
\scriptfont2=\sevensy%
\scriptscriptfont0=\sixrm%
\scriptscriptfont1=\sixi%
\scriptscriptfont2=\sixsy%
\def\it{\fam\itfam\tenit}%
\def\tt{\fam\ttfam\tentt}%
\def\sl{\fam\slfam\tensl}%
\scriptfont\bffam=\sevenbf%
\scriptscriptfont\bffam=\sixbf%
\def\bf{\fam\bffam\tenbf}%
\normalbaselineskip=18pt%
\normalbaselines\rm}%

\font\ninerm=cmr9%
\font\ninebf=cmbx9%
\font\nineit=cmti9%
\font\ninesy=cmsy9%
\font\ninei=cmmi9%
\font\ninett=cmtt9%
\font\ninesl=cmsl9%

\def\ninepoint{\def\rm{\fam0\ninerm}%
\textfont0=\ninerm%
\textfont1=\ninei%
\textfont2=\ninesy%
\textfont\itfam=\nineit%
\textfont\slfam=\ninesl%
\textfont\ttfam=\ninett%
\textfont\bffam=\ninebf%
\scriptfont0=\sixrm%
\scriptfont1=\sixi%
\scriptfont2=\sixsy%
\def\it{\fam\itfam\nineit}%
\def\tt{\fam\ttfam\ninett}%
\def\sl{\fam\slfam\ninesl}%
\scriptfont\bffam=\sixbf%
\scriptscriptfont\bffam=\fivebf%
\def\bf{\fam\bffam\ninebf}%
\normalbaselineskip=16pt%
\normalbaselines\rm}%

\font\eightrm=cmr8%
\font\eighti=cmmi8%
\font\eightsy=cmsy8%
\font\eightbf=cmbx8%
\font\eighttt=cmtt8%
\font\eightit=cmti8%
\font\eightsl=cmsl8%

\def\eightpoint{\def\rm{\fam0\eightrm}%
\textfont0=\eightrm%
\textfont1=\eighti%
\textfont2=\eightsy%
\textfont\itfam=\eightit%
\textfont\slfam=\eightsl%
\textfont\ttfam=\eighttt%
\textfont\bffam=\eightbf%
\scriptfont0=\sixrm%
\scriptfont1=\sixi%
\scriptfont2=\sixsy%
\scriptscriptfont0=\fiverm%
\scriptscriptfont1=\fivei%
\scriptscriptfont2=\fivesy%
\def\it{\fam\itfam\eightit}%
\def\tt{\fam\ttfam\eighttt}%
\def\sl{\fam\slfam\eightsl}%
\scriptscriptfont\bffam=\fivebf%
\def\bf{\fam\bffam\eightbf}%
\normalbaselineskip=14pt%
\normalbaselines\rm}%

\font\sevenrm=cmr7%
\font\seveni=cmmi7%
\font\sevensy=cmsy7%
\font\sevenbf=cmbx7%

\def\sevenpoint{%
   \def\rm{\sevenrm}\def\bf{\sevenbf}%
   \def\smc{\sevensmc}\baselineskip=12pt\rm}%

\font\sixrm=cmr6%
\font\sixi=cmmi6%
\font\sixsy=cmsy6%
\font\sixbf=cmbx6%

\fontdimen13\tensy=2.6pt%
\fontdimen14\tensy=2.6pt%
\fontdimen15\tensy=2.6pt%
\fontdimen16\tensy=1.2pt%
\fontdimen17\tensy=1.2pt%
\fontdimen18\tensy=1.2pt%

\def\tenf{\tenpoint}%
\def\ninef{\ninepoint}%
%



\def\section#1{\goodbreak\vskip 3pc plus 6pt minus 3pt\leftskip=-2pc
   \global\advance\sectnum by 1\eqnumber=1\subsectnum=0%
\global\examplnumber=1\figrnumber=1\propnumber=1\defnumber=1\lemmanumber=1\assumptionnumber=1 \conditionnumber =1%
   \line{\hfuzz=1pc{\hbox to 3pc{\bf 
   \vtop{\hfuzz=1pc\hsize=38pc\hyphenpenalty=10000\noindent\uppercase{\the\sectnum.\quad #1}}\hss}}
			\hfill}
			\leftskip=0pc\nobreak\tenf
			\vskip 1pc plus 4pt minus 2pt\noindent\ignorespaces}
\def\subsection#1{\noindent\leftskip=0pc\tenf
   \goodbreak\vskip 1pc plus 4pt minus 2pt
               \global\advance\subsectnum by 1
   \line{\hfuzz=1pc{\hbox to 3pc{\bf \the\sectnum.\the\subsectnum.
   \vtop{\hfuzz=1pc\hsize=38pc\hyphenpenalty=10000\noindent{\bf #1}}\hss}}
                        \hfill}
   \leftskip=0pc\nobreak\tenf
                        \vskip 1pc plus 4pt minus 2pt\nobreak\noindent\ignorespaces}



\def\texshopbox#1{\boxtext{462pt}{\vskip-1.5pc\pshade{\vskip-1.0pc#1\vskip-2.0pc}}}


\input miniltx

\ifx\pdfoutput\undefined
  \def\Gin@driver{dvips.def} 
\else
  \def\Gin@driver{pdftex.def} 
\fi

\input graphicx.sty
\resetatcatcode

\long\def\fig#1#2#3{\vbox{\vskip1pc\vskip#1
\prevdepth=12pt \baselineskip=12pt
\vskip1pc
\hbox to\hsize{\hfill\vtop{\hsize=30pc\noindent{\eightbf Figure #2\ }
{\eightpoint#3}}\hfill}}}

\def\show#1{}

\rightheadline{\botmark}

\pageno=1

\def\longpapertitle#1#2#3{{\bf \centerline{\helbigb
{#1}}}\medskip{\bf \centerline{\helbigb
{#2}}}\bigskip{\bf \centerline{
{#3}}}\bigskip}

\vskip-3pc

\def\xstar{X^{\raise0.04pt\hbox{\sevenpoint *}} }

\def\jstar{J^{\raise0.04pt\hbox{\sevenpoint *}} }
\def\qstar{Q^{\raise0.04pt\hbox{\sevenpoint *}} }

\rightheadline{\botmark}

\pageno=1

\rightheadline{\botmark}

\pageno=1

\rightheadline{\botmark}

\pn {\bf March 2017 (Revised May 2017)}\hfill{\bf  Report LIDS-P-3506}%
\bigskip \bigskip\bigskip

\bigskip

\def\longpapertitle#1#2#3{{\bf \centerline{\helbigb
{#1}}}\medskip{\bf \centerline{\helbigb
{#2}}}\bigskip{\bf \centerline{
{#3}}}\bigskip}

\vskip-3pc

\longpapertitle{Stable  Optimal Control and}{Semicontractive Dynamic Programming}{ {Dimitri P.\ Bertsekas\footnote{\dag}{\ninepoint  Dimitri Bertsekas is with the Dept.\ of Electr.\ Engineering and
Comp.\ Science, and the Laboratory for Information and Decision Systems, M.I.T., Cambridge, Mass., 02139.}  }}


\centerline{\bf Abstract}

We consider discrete-time infinite horizon deterministic optimal control problems with nonnegative cost per stage, and a destination that is cost-free and absorbing. The classical linear-quadratic regulator problem is a special case. Our assumptions are very general, and allow the possibility that the optimal policy may not be stabilizing the system, e.g., may not reach the destination either asymptotically or in a finite number of steps. We introduce a new unifying notion of stable feedback policy, based on perturbation of the cost per stage, which in addition to implying convergence of the generated states to the destination, quantifies the speed of convergence. We consider the properties of two distinct cost functions: $\jstar$, the overall optimal, and $\hat J$, the restricted optimal over just the stable policies. Different classes of stable policies (with different speeds of convergence) may yield different values of $\hat J$. We show that for any class of stable policies, $\hat J$ is a solution of Bellman's equation, and we characterize the smallest and the largest solutions: they are $\jstar$, and $J^+$, the restricted optimal cost function over the class of (finitely) terminating policies. We also characterize the regions of convergence of various modified versions of value and policy iteration algorithms, as substitutes for the standard algorithms, which may not work in general.

\vskip-1.5pc

\section{Introduction}

\vskip-0.5pc

\pn In this paper we consider a deterministic discrete-time infinite horizon optimal control problem involving the system
$$x_{k+1}=f(x_k,u_k),\qquad k=0,1,\ldots,\xdef\docsys{\lab}\eqnum\show{oneo}$$
where $x_k$ and $u_k$ are the state and control at stage $k$, which belong to sets $X$ and $U$, referred to as the state and control spaces, respectively, and $f:X\times U\mapsto X$ is a given function. The control $u_k$ must be chosen from a nonempty constraint set $U(x_k)\subset U$ that may depend on the current state $x_k$. The cost for the $k$th stage, $g(x_k,u_k)$, is assumed nonnegative and possibly extended  real-valued:
$$0\le g(x_k,u_k)\le\infty,\qquad \forall\ x_k\in X,\ u_k\in U(x_k),\ k=0,1,\ldots.\xdef\nonnegcost{\lab}\eqnum\show{oneo}$$
A cost per stage that is extended real-valued may be useful in modeling conveniently additional state and control constraints. We assume that $X$ contains a special state, denoted $t$, which is referred to as the {\it destination\/}, and is cost-free and absorbing:
$$f(t,u)=t,\qquad g(t,u)=0,\qquad \forall\ u\in U(t).\eqnum\show{oneo}$$

Our terminology aims to emphasize the connection with classical problems of control where $X$ and $U$ are the finite-dimensional Euclidean spaces $X=\rn$, $U=\re^m$, and the destination is identified with the origin of $\rn$. There the essence of the problem is to reach or asymptotically approach the origin at minimum cost. A  special case is the classical infinite horizon linear-quadratic regulator problem. However, our formulation also includes shortest path problems with continuous as well as discrete spaces; for example the classical shortest path problem, where $X$ consists of the nodes of a directed graph, and the problem is to reach the destination from every other node with a minimum length path.

We are interested in  feedback policies of the form $\p=\{\m_0,\m_1,\ldots\}$, where each $\m_k$ is a function mapping $x\in X $ into the control $\m_k(x)\in U(x)$. The set of all  policies is denoted by $\P$. Policies of the form $\p=\{\m,\m,\ldots\}$ are called {\it stationary\/}, and will be denoted by $\m$, when confusion cannot arise. 

Given an initial state $x_0$, a policy $\p=\{\m_0,\m_1,\ldots\}$ when applied to the system \docsys, generates a unique sequence of state-control pairs $\big(x_k,\m_k(x_k)\big)$, $k=0,1,\ldots,$ with cost 
$$J_\p(x_0)= \sum_{k=0}^{\infty}
g\bl(x_k,\mu_k(x_k)\br),\qquad x_0\in X ,$$
[the series converges to some number in $[0,\infty]$ thanks to the nonnegativity assumption \nonnegcost]. We view $J_\p$ as a function over $X $, and we refer to it as the cost function of $\p$. For a stationary policy $\m$, the corresponding cost function is denoted by $J_\m$.
The optimal cost function is defined as
$$\jstar(x)=\inf_{\p\in\P}J_\p(x),\qquad x\in X ,$$
and a policy $\p^*$ is said to be optimal if 
$J_{\p^*}(x)=\jstar(x)$ for all $x\in X .$ The optimal cost $\jstar(x)$ is identical to the optimal cost attained when  starting at $x$ and using open-loop sequences $\{u_0,u_1,\ldots\}$, but in this paper we consider the broader class of feedback policies to be consistent with the formalism and analysis of infinite horizon DP.

We denote by ${\cal E}^+(X)$ the set of functions $J:X\mapsto[0,\infty]$. All equations, inequalities, limit and minimization operations involving functions from this set are meant to be pointwise. In our analysis, we will use the set of functions
$${\cal J}=\big\{J\in{\cal E}^+(X) \mid J(t)=0\big\}.$$
Since $t$ is cost-free and absorbing, this set contains $J_\p$ of all $\p\in \P$, as well as $\jstar$.  

It is well known that under the cost nonnegativity assumption \nonnegcost, $\jstar$ satisfies Bellman's equation:
$$\jstar(x)=\inf_{u\in U(x)}\big\{g(x,u)+\jstar\big(f(x,u)\big)\big\},\qquad x\in X ,$$
and that an optimal stationary policy (if it exists) may be obtained through the minimization in the right side of this equation (cf.\ Prop.\ 2.1 in the next section). One also hopes to obtain $\jstar$ by means of value iteration (VI for short), which starting from some function $J_0\in {\cal J}$, generates a sequence of functions $\{J_k\}\subset  {\cal J}$ according to 
$$J_{k+1}(x)=\inf_{u\in U(x)}\big\{g(x,u)+J_k\big(f(x,u)\big)\big\},\qquad x\in X ,\ k=0,1,\ldots.\xdef\vieq{\lab}\eqnum\show{oneo}$$
However, $\{J_k\}$ may not always converge to $\jstar$ because, among other reasons, Bellman's equation may have multiple solutions within ${\cal J}$.

Another possibility to obtain $\jstar$ and an optimal policy is through policy iteration (PI for short), which starting from a stationary policy $\m^0$, generates a sequence of stationary policies $\{\m^k\}$ via a sequence of policy evaluations to obtain $J_{\m^k}$ from the equation
$$J_{\m^k}(x)=g\big(x,\m^k(x)\big)+J_{\m^k}\big(f\big(x,\m^k(x)\big)\big),\qquad x\in X ,\xdef\poleval{\lab}\eqnum\show{oneo}$$
interleaved with policy improvements to obtain $\m^{k+1}$ from $J_{\m^k}$ according to
$$\m^{k+1}(x)\in \argmin_{u\in U(x)}\big\{g(x,u)+J_{\m^k}\big(f(x,u)\big)\big\},\qquad x\in X .\xdef
\polimprove{\lab}\eqnum\show{oneo}$$
We note that $J_{\m^k}$ satisfies Eq.\ \poleval\ (cf.\ Prop.\ 2.1 in the next section). Moreover, when referring to PI, we  assume that the minimum in Eq.\ \polimprove\ is attained for all $x\in X $, which is true under some compactness condition on  $U(x)$ or the level sets of the function $g(x,\cdot)+J_k\big(f(x,\cdot)\big)$, or both (see e.g., [Ber12], Ch.\ 4).

The uniqueness of solution of Bellman's equation within ${\cal J}$, and the convergence of VI to $\jstar$ and of PI to an optimal policy have been investigated in a recent paper by the author [Ber15a] under conditions guaranteeing that $\jstar$ is the unique solution of Bellman's equation within ${\cal J}$.  This paper also gives many references from the field of adaptive dynamic programming, where the continuous-spaces version of our problem is often used as the starting point for analysis and algorithmic development; see e.g., the book [VVL13], the papers [JiJ14], [LiW13], the survey
papers in the edited volumes [SBP04] and [LeL13], and the special issue  [LLL08]. Our purpose here is to consider the problem under weaker conditions and to make the connection with notions of stability. This is a more complicated case, where Bellman's equation need not have a unique solution within ${\cal J}$, while the VI and PI algorithms may be unreliable. However, several of the favorable results of [Ber15a] will be obtained as special cases of the results of this paper; see Section 3. The type of behavior that we are trying to quantify is described in the following example.\footnote{\dag}{\ninepoint  In this example and later, our standard notational convention is that all vectors in $\rn$ are viewed as column vectors. The real line is denoted by $\re$. A prime denotes transposition, so inner product of two vectors $x$ and $y$ is defined by $x'y$, and the norm is $\|x\|=\sqrt{x'x}$.}

\xdef\examplelq{\exampl}\examplnum\show{examplo}

\xdef \figlinquad{\figr}\figrnum\show{myfigure}
 
\topinsert
\centerline{\hskip0pc\includegraphics[width=3.2in]{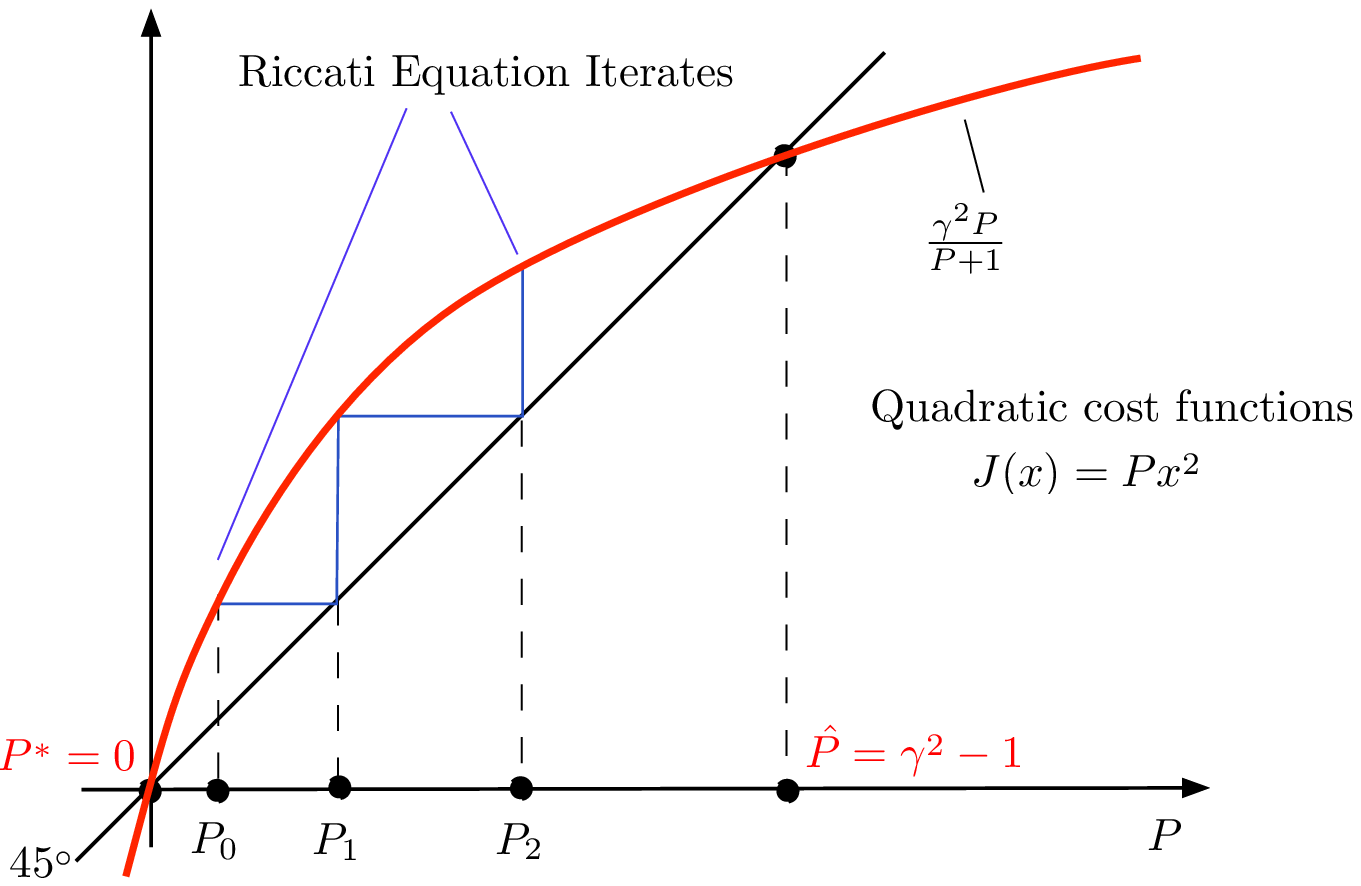}}
\fig{-1.2pc}{\figlinquad} {Illustration of the  behavior of the Riccati equation for the linear-quadratic problem of Example \examplelq, where the detectability assumption is not satisfied. The solutions of the Riccati equation are $P=0$ (corresponds to the optimal cost) and $\skew4\hat P=\g^2-1$ (corresponds to the optimal cost that can be achieved with linear stable control laws).}
\endinsert

 \beginexample{\examplelq\ (Linear-Quadratic Problem)}Consider a linear system and a quadratic cost:
$$x_{k+1}=Ax_k+Bu_k,\qquad g(x_k,u_k)=x_k'Qx_k+u_k'Ru_k,$$
where $X=\rn$, $U=\re^m$, $A$, $B$, $Q$, and $R$ are given matrices, with $Q$ being positive semidefinite symmetric and $R$ being positive definite symmetric. The classical results for this problem assume that:
\nitem{(a)} The pair $(A,B)$ is stabilizable [i.e., there exists a linear policy $\m(x)=Lx$ such that the closed-loop system $x_{k+1}=(A+BL)x_k$ is asymptotically stable].
\nitem{(b)} The pair $(A,C)$, where $Q=C'C$, is detectable (i.e., if $u_k\to0$ and $Cx_k\to0$ then it follows that $x_k\to0$).
\smskip
\pn Under these assumptions, it is well-known (see e.g., optimal control and estimation textbooks such as Anderson and Moore [AnM07], or the author's [Ber17], Section 3.1) that $J^*$ has the form $J^*(x)=x'Px$, where $P$ is the unique positive semidefinite solution of the algebraic Riccati equation
$$P = A'\bl( P-PB(B'PB + R)^{-1}B'P\br) A+Q.\xdef\riccati{\lab}\eqnum\show{threeo}$$
Furthermore the VI algorithm converges to $J^*$ starting from any positive semidefinite quadratic function. Moreover the PI algorithm, starting from a linear stable policy, yields $J^*$ and a linear stable optimal policy in the limit as first shown by Kleinman [Kle68].

To see what may happen when the preceding detectability condition is not satisfied, consider the scalar  system 
$$x_{k+1}=\g x_k+u_k,$$
where $\g>1$ and the  cost per stage is $g(x,u)=u^2$. 
Here we have $J^*(x)\equiv0$, while the policy $\m^*(x)\equiv0$ is optimal. This policy is not stable (for any sensible definition of stability), which is not inconsistent with optimality, since nonzero states are not penalized in this problem. The  algebraic Riccati equation \riccati\ for this case is
$$P={\g^2 P\over P+1},$$
and has {\it two nonnegative solutions\/}: $P^*=0$ and $\skew4\hat P=\g^2-1$ (see Fig.\ \figlinquad).  It turns out that $\skew4\hat P$ has an interesting interpretation: $\skew5\hat J(x)=\skew4\hat P x^2$ is the optimal cost that can be achieved using a linear stable control law, starting from $x$ (see the analysis of [Ber17], Example 3.1.1). Moreover the VI algorithm, which generates the sequence $J_k(x)=P_kx^2$, with $P_k$ obtained by the Riccati equation iteration 
$$P_{k+1}={\g^2 P_k\over P_k+1},$$ 
converges to $\skew5\hat J$ when started  from any $P_0>0$, and stays at the zero function $J^*$ when started  from $P_0=0$ (see Fig.\ \figlinquad). Another interesting fact is that the PI algorithm, when started from a linear stable policy, yields in the limit $\skew5\hat J$ (not $J^*$) and the policy that is optimal within the class of linear stable policies (which turns out to be $\hat \m(x)={1-\g^2\over\g}x$); see [Ber17], Section 3.1, for a verification, and Example 3.1 for an analysis of the multidimensional case.\footnote{\dag}{\ninepoint As an example of what may happen without stabilizability,  consider the case when the system is instead $x_{k+1}=\g x_k$. Then the Riccati equation becomes $P=\g^2 P$ and has $P^*=0$ as its unique solution. However, the Riccati equation iteration $P_{k+1}=\g^2 P_k$ diverges to $\infty$ starting from any $P_0>0$. Also, qualitatively similar behavior is obtained when there is a discount factor $\a\in(0,1)$. The Riccati equation takes the form
$$P = A'\bl( \a P-\a^2 PB(\a B'PB + R)^{-1}B'P\br) A+Q,$$
 and for the given system and cost per stage, it has two solutions, $P^*=0$ and $\skew4\hat P={\a\g^2-1\over \a}$. The VI algorithm converges to $\skew4\hat P$ starting from any $P>0$. While the line of analysis of the present paper does not apply to discounted problems, a related analysis is given in the paper [Ber15b], using the idea of a regular policy.} 

We note that the set of solutions of the Riccati equation has been extensively investigated starting with the papers by Willems [Wil71] and Kucera [Kuc72], [Kuc73], which were followed up by several other works (see the book by Lancaster and Rodman [LaR95] for a comprehensive treatment).  In these works, the ``largest" solution of the Riccati equation is referred to as the ``stabilizing" solution, and the stability of the corresponding policy is shown, although the author could not find an explicit statement regarding the optimality of this policy within the class of all linear stable policies. Also the lines of analysis of these works are tied to the structure of the linear-quadratic problem and are unrelated to the analysis of the present paper.
\endexample

\vskip-1pc

There are also other interesting deterministic optimal control examples where Bellman's equation, and the VI and PI algorithms exhibit unusual behavior, including several types of shortest path problems (see e.g., [Ber14], [Ber15a], [BeY16], and the subsequent Example 4.1).
This is typical of {\it semicontractive DP} theory, which is a central focal point of the author's abstract DP monograph [Ber13], and followup work [Ber15b]. The present paper is inspired by the analytical methods of this theory. In semicontractive models, roughly speaking, policies are divided into those that are ``regular" in the sense that they are ``well-behaved" with respect to the VI algorithm, and those that are ``irregular." The optimal cost function over the ``regular" policies (under suitable conditions) is a solution of Bellman's equation, and can be found by the VI and PI algorithm, even under conditions where these algorithms may fail to find the optimal cost function $\jstar$. Regularity in the sense of semicontractive DP corresponds to stability in the specialized context of deterministic optimal control considered here.

In this paper we address the phenomena illustrated by the  linear-quadratic Example \examplelq\ in the more general setting where the system may be nonlinear and the cost function may be nonquadratic. Our method of analysis is to introduce a cost perturbation that involves a penalty for excursions of the state from the destination, thus resulting in a better-behaved problem. The type of perturbation used determines in turn the class of stable policies. A key aspect of our definition of a stable policy  (as given in the next section) is that in addition to implying convergence of the generated states to the destination, it quantifies the speed of convergence. 

A simpler approach, which involves perturbation by a constant function, has been used in the monograph [Ber13], and also in the paper by Bertsekas and Yu [BeY16]. The latter paper analyzes similarly unusual behavior in finite-state finite-control stochastic shortest path problems, where the cost per stage can take both positive and negative values (for such problems the anomalies are even more acute,  including the possibility that $\jstar$ may not solve Bellman's equation). 

In the analysis of the present paper, the optimal policies of the perturbed problem are stable policies, and in the limit as the perturbation diminishes to 0, the corresponding optimal cost function converges to $\hat J$, the optimal cost function over stable policies (not to $\jstar$). Our central result is that $\hat J$ is the unique solution of Bellman's equation within a set of functions in ${\cal J}$ that majorize $\hat J$. Moreover, the VI algorithm converges to $\hat J$ when started from within this set. In addition, if $J^+$,  the optimal cost function over the class of (finitely) terminating policies belongs to ${\cal J}$, then $J^+$  is the largest solution of Bellman's equation within ${\cal J}$. These facts are shown in Section 3, including a treatment of the multidimensional version of the linear-quadratic problem of Example \examplelq. In Section 3, we also consider the favorable special case where $\jstar=J^+$, and we develop the convergence properties of VI for this case. In Section 4 we consider PI algorithms, including a perturbed version (Section 4.2).

\vskip-1.5pc

\section{Stable Policies}

\pn In this section, we will lay the groundwork for our analysis and introduce the notion of a stable policy. To this end, we will use some classical results for optimal control with nonnegative cost per stage, which stem from the original work of Strauch [Str66]. For textbook accounts we refer to [BeS78], [Put94], [Ber12], and for a more abstract development, we refer to the monograph [Ber13].

The following proposition gives the results that we will need (see [BeS78], Props.\ 5.2, 5.4, and 5.10, [Ber12], Props.\ 4.1.1, 4.1.3, 4.1.5, 4.1.9, or [Ber13], Props.\ 4.3.3, 4.3.9, and 4.3.14). Actually, these results hold for stochastic infinite horizon DP problems with nonnegative cost per stage, and do not depend on the favorable structure of this paper (a deterministic problem with a cost-free and absorbing destination).


\xdef\propnegdp{\propn}\propnum\show{myproposition}

\texshopbox{\proposition{\propnegdp:}The following hold:
\nitem{(a)} $\jstar$ is a solution of Bellman's equation and  if $J\in{\cal E}^+(X)$ is another solution, i.e., $J$ satisfies
$$J(x)= \inf_{u\in U(x)}\big\{g(x,u)+J\big(f(x,u)\big)\big\},\qquad \forall\ x\in X ,\xdef\bellmaneq{\lab}\eqnum\show{oneo}$$
 then  $\jstar\le J$.
\nitem{(b)} For all stationary policies $\m$, $J_\m$ is a solution of the equation 
$$J(x)=g\big(x,\m(x)\big)+J\big(f\big(x,\m(x)\big)\big),\qquad \forall\ x\in X,$$ 
and if $J\in{\cal E}^+(X)$ is another solution, then $J_\m\le J$.
\nitem{(c)} A stationary policy $\m^*$ is optimal if and only if 
$$\m^*(x)\in \argmin_{u\in U(x)}\big\{g(x,u)+\jstar\big(f(x,u)\big)\big\},\qquad \forall\ x\in X.$$
\nitem{(d)} Let $\{\bar J_k\}$ be the sequence generated by the VI  algorithm \vieq\ starting from the zero function $\bar J_0(x)\equiv0$. If $U$ is a metric space and the sets 
$$U_k(x,\l) =\big\{ u\in U(x)\big|\ g(x,u)+\bar J_k\big(f(x,u)\big)\le\l\big\}\eqnum\show{oneo}$$
are compact for all $x\in X$, $\l\in \re$, and $k\ge0$,  then there exists at least one optimal stationary policy, and we have $J_k\to\jstar$ for every sequence generated by VI starting from a function $J_0\in{\cal E}^+(X)$ with $J_0\le \jstar$.
\nitem{(e)} For every $\e>0$ there exists an $\e$-optimal  policy, i.e., a policy $\p_{\e}$ such that
$$J_{\p_\e}(x)\le \jstar(x)+\e,\qquad\forall\ x\in X.$$
}

\old{
DO WE NEED? The following proposition depends not only on the nonnegativity of the cost per stage, but also on the deterministic character of the problem.
\xdef\proppolchar{\propn}\propnum\show{myproposition}
\texshopbox{\proposition{\proppolchar:}Let $\p=\{\m_0,\m_1,\ldots\}$ be a policy, and for a given initial state $x_0\in X$, let $\{x_k\}$ be the sequence of  states generated by starting from $x_0$ and using $\p$.
\nitem{(a)} If $J_\p(x_0)<\infty$, then $J_{\p_k}(x_k)\downarrow0$, where $\p_k$ is the  policy $\{\m_k,\m_{k+1},\ldots\}$.
\nitem{(b)} If $\p$ is  stationary of the form $\{\m,\m,\ldots\}$ and $J_\m(x_0)<\infty$, then $J_\m(x_k)\downarrow0$. 
}
\proof (a) We have by definition
$$J_{\p_m}(x_m)=g\big(x_m,\m_m(x_m)\big)+J_{\p_{m+1}}(x_{m+1}),\qquad m=0,1,\ldots.$$
By repeatedly applying this relation, we obtain
$$J_\p(x_0)=\sum_{m=0}^{k-1}g\big(x_m,\m_m(x_m)\big)+J_{\p_k}(x_{k}),\qquad k=0,1,\ldots.$$
Since $g$ is nonnegative, $\sum_{m=0}^{k-1}g\big(x_m,\m_m(x_m)\big)$ is monotonically nondecreasing, and it follows that $J_{\p_k}(x_k)$ is monotonically nonincreasing and real valued since $J_\p(x_0)<\infty$.
Moreover, by taking the limit as $k\to\infty$ we obtain $\lim_{k\to\infty}J_{\p_k}(x_{k})=0$.
\smskip
\pn (b) This is a special case of part (a), with  $\p=\{\m,\m,\ldots\}$. \qed
}

We introduce a {\it forcing function} $p:X\mapsto[0,\infty)$ such that 
$$p(t)=0,\qquad p(x)>0,\qquad \forall\ x\ne t,$$ 
and the {\it $p$-$\d$-perturbed optimal control problem\/}, where $\d>0$ is a given scalar. This is the same problem as the original, except that the cost per stage is changed to
$$g(x,u)+\d p(x).$$
We denote by $J_{\p,p,\d}$ the cost function of a policy $\p\in\P$ in the $p$-$\d$-perturbed problem: $$J_{\p,p,\d}(x_0)=J_\p(x_0)+\d\sum_{k=0}^\infty  p(x_k),\xdef\costexpr{\lab}\eqnum\show{oneo}
$$
where $\{x_k\}$ is the sequence generated starting from $x_0$ and using $\p$. We also denote by $\hat J_{p,\d}$, the corresponding optimal cost function, $\hat J_{p,\d}=\inf_{\p\in\Pi}J_{\p,p,\d}$.
We  introduce a notion of stability involving the $p$-$\d$-perturbed problem.

\xdef\definitionstable{\defn}\defnum\show{myproposition}

\texshopbox{\definition{\definitionstable:}Let $p$ be a given forcing function. For a state $x\in X$, we say that a policy $\p$ is {\it $p$-stable from $x$} if for all $\d>0$ we have 
$$J_{\p,p,\d}(x)<\infty.$$
The set of all such policies is denoted by $\P_{p,x}$. We define the {\it restricted optimal cost function} over $\P_{p,x}$ by
$$\hat J_p(x)=\inf_{\p\in\P_{p,x}}J_\p(x),\qquad x\in X.\xdef\jhatpjdef{\lab}\eqnum\show{oneo}$$
We say that $\p$ is {\it $p$-stable} if $\p\in \P_{p,x}$ simultaneously from all $x\in X$ such that $\P_{p,x}\ne\emptyset$. The set of all $p$-stable policies is denoted by $\P_p$.
}

The preceding definition of a $p$-stable policy is novel within the general context of this paper, and is inspired from a notion of regularity, which is central in the theory of semicontractive DP; see [Ber13] and the related subsequent papers [Ber14], [Ber15b], [Ber16]. Note that
the set $\P_{p,x}$ depends on the forcing function $p$. As an example, let $X=\rn$ and 
$$p(x)=\|x\|^\r,$$ where $\r>0$ is a scalar.  Then roughly speaking, $\r$ quantifies the rate at which the destination is approached using the $p$-stable policies. In particular, the policies $\p\in \P_{p,x_0}$ are the ones that force $x_k$ towards $0$ at a rate faster than $O(1/k^\r)$, so slower policies would be excluded from $\P_{p,x_0}$.

Let us make some observations regarding $p$-stability:

\nitem{(a)} {\it Equivalent definition of $p$-stability\/}: Given any policy $\p$ and state $x_0\in X$, from Eq.\ \costexpr\ 
 it follows that
$$\p\in \P_{p,x_0}\qquad \hbox{if and only if}\qquad J_\p(x_0)<\infty\hbox{ and }\sum_{k=0}^\infty p(x_k)<\infty,\xdef\equivstab{\lab}\eqnum\show{oneo}$$
where $\{x_k\}$ is the sequence generated starting from $x_0$ and using $\p$. Since the right-hand side of the preceding relation does not depend on $\d$, it also follows that an equivalent definition of a policy $\p$ that is $p$-stable from $x$ is that $J_{\p,p,\d}(x)<\infty$ for some $\d>0$ (rather than all $\d>0$).

\nitem{(b)}  {\it Approximation property of  $J_{\p,p,\d}(x)$\/}: Consider a pair $(\p,x_0)$ with $\p\in \P_{p,x_0}$. By taking the limit as $\d\downarrow0$ in the expression
$$J_{\p,p,\d}(x_0)=J_\p(x_0)+\d\sum_{k=0}^\infty  p(x_k),$$
[cf.\ Eq.\ \costexpr] and by using Eq.\ \equivstab, it follows that
$$\lim_{\d\downarrow0}J_{\p,p,\d}(x_0)=J_{\p}(x_0),\qquad  \forall\ \hbox{pairs $(\p,x_0)$ with $\p\in \P_{p,x_0}$}.\xdef\limdeltazero{\lab}\eqnum\show{oneo}$$
From this equation, we have that if $\p\in \P_{p,x}$, then $J_{\p,p,\d}(x)$ is finite and differs from $J_\p(x)$ by $O(\d)$. By contrast, if $\p\notin \P_{p,x}$, then $J_{\p,p,\d}(x)=\infty$ by the definition of $p$-stability. 

\nitem{(c)}  {\it Limiting property of  $\hat J_p(x_k)$\/}: Consider a pair $(\p,x_0)$ with $\p\in \P_{p,x_0}$. By breaking down $J_{\p,p,\d}(x_0)$ into the sum of the costs of the first $k$ stages and the remaining stages, we have
$$J_{\p,p,\d}(x_0)=\sum_{m=0}^{k-1}  g\big(x_m,\m_m(x_m)\big)+\d\sum_{m=0}^{k-1}  p(x_m)+J_{\p_k,p,\d}(x_k),\qquad  \forall\ \d>0,\ k>0,\eqnum\show{oneo}$$
where $\{x_k\}$ is the sequence generated starting from $x_0$ and using $\p$, and $\p_k$ is the policy $\{\m_k,\m_{k+1},\ldots\}$. By taking the limit as $k\to\infty$ and using Eq.\ \costexpr, it follows that
$$\lim_{k\to\infty}J_{\p_k,p,\d}(x_k)= 0,\qquad \forall\ \hbox{pairs $(\p,x_0)$ with $\p\in \P_{p,x_0}$},\ \d>0.\eqnum\show{oneo}$$
Also, since $\hat J_p(x_k)\le \hat J_{p,\d}(x_k)\le J_{\p_k,p,\d}(x_k)$, it follows that
$$\lim_{k\to\infty}J_{p,\d}(x_k)=0,\qquad \lim_{k\to\infty}\hat J_p(x_k)=0,\qquad \hbox{$\forall\ (\p,x_0)$ with $x_0\in X$ and $\p\in\P_{p,x_0}$, $\d>0$}.\xdef\limprophatj{\lab}\eqnum\show{oneo}$$

\subsubsection{Terminating Policies and Controllability}

\pn An important special case is when $p$ is equal to the function
$$p^+(x)=\cases{0&if $x=t$,\cr
1&if $x\ne t$.\cr}\xdef\terminatingp{\lab}\eqnum\show{oneo}$$
For $p=p^+$, a policy $\p$ is $p^+$-stable from $x$ if and only if it is {\it terminating from $x$\/}, i.e., reaches $t$ in a finite number of steps starting from $x$ [cf.\ Eq.\ \equivstab]. 
The set of terminating policies from $x$ is denoted by ${\P}_x^+$ and it is contained within every other set of $p$-stable policies $\P_{p,x}$, as can be seen from Eq.\ \equivstab. As a result, the restricted optimal cost function over ${\P}_x^+$, 
$$J^+(x)=\inf_{\p\in{\P}_x^+}J_\p(x),\qquad x\in X,$$
satisfies
$\jstar(x)\le \hat J_p(x)\le J^+(x)$ for all $x\in X.$
A policy $\p$ is said to be {\it terminating} if it is simultaneously terminating from all $x\in X$ such that ${\P}_x^+\ne \emptyset$. The set of all terminating policies is denoted by $\P^+$.

Note that if the state space $X$ is finite, we have for every forcing function $p$
$$\ubar \b\, p^+(x)\le p(x)\le \bar \b\, p^+(x),\qquad\forall\ x\in X,$$
for some scalars $\ubar \b,\bar \b>0$. As a result it can be seen that $\P_{p,x}=\P^+_x$ and $\hat J_p=J^+$, so in effect the case where $p=p^+$ is the only case of interest for finite-state problems.

The notion of a terminating policy is related to the  notion of {\it controllability\/}. In classical control theory terms, the system $x_{k+1}=f(x_k,u_k)$ is said to be completely controllable if for every $x_0\in X$, there exists a policy that drives the state $x_k$ to the destination in a finite number of steps. This notion of controllability is equivalent to  the existence of a terminating policy from each $x\in X$.

One of our main results, to be shown in the next section, is that $\jstar$, $\hat J_p$, and $J^+$ are solutions of Bellman's equation,  with $\jstar$ being the ``smallest" solution and $J^+$ being the ``largest" solution within ${\cal J}$.
The most favorable situation arises when $\jstar= J^+$, in which case $\jstar$ is the unique solution of Bellman's equation within ${\cal J}$. Moreover, in this case it will be shown that the VI algorithm converges to $\jstar$ starting with any $J_0\in {\cal J}$ with $J_0\ge \jstar$ (see the subsequent Prop.\ 3.5), and the PI algorithm converges to $\jstar$ as well (see Section 4.1). This special case has been discussed in the paper [Ber15a].

\vskip-1.5pc

\section{Restricted Optimization over Stable Policies}

\pn For a given forcing function $p$, we denote by 
$\widehat X_p$ the effective domain of $\hat J_p$, the set of all $x$ where $\hat J_p$ is finite,
$$\widehat X_p=\big\{x\in X\mid \hat J_p(x)<\infty\big\}.$$
Since $\hat J_p(x)<\infty$ if and only if $\P_{p,x}\ne\emptyset$ [cf.\ Eqs.\ \jhatpjdef\ and \equivstab], or equivalently $J_{\p,p,\d}(x)<\infty$ for some $\p$ and all $\d>0$, 
it follows that $\widehat X_p$ is also the effective domain of $\hat J_{p,\d}$, 
$$\widehat X_p=\big\{x\in X\mid \P_{p,x}\ne\emptyset\}=\big\{x\in X\mid \hat J_{p,\d}(x)<\infty\big\},\qquad \forall\ \d>0.$$
Note that $\widehat X_p$ may depend on $p$ and may be a strict subset of the effective domain of $\jstar$, which is denoted by
$$\xstar=\big\{x\in X \mid \jstar(x)<\infty\big\}.$$
The reason is that there may exist a policy $\p$ such that $J_\p(x)<\infty$, even when there is no $p$-stable policy from $x$.

Our first objective is to show that as $\d\downarrow0$, the $p$-$\d$-perturbed optimal cost function $\hat J_{p,\d}$ converges to the restricted optimal cost function $\hat J_p$.

\xdef\propdeltaper{\propn}\propnum\show{myproposition}

\texshopbox{\proposition{\propdeltaper\ (Approximation Property of $\hat J_{p,\d}$):}Let $p$ be a given forcing function and $\d>0$.
\nitem{(a)} We have
$$J_{\p,p,\d}(x)= J_\p(x)+w_{\p,p,\d}(x),\qquad \forall\ x\in X,\ \p\in \P_{p,x},\xdef\approxform{\lab}\eqnum\show{oneo}$$
where $w_{\p,p,\d}$ is a function such that $\lim_{\d\downarrow0}w_{\p,p,\d}(x)=0$ for all $x\in X$.
\nitem{(b)}We have 
$$\lim_{\d\downarrow0}\hat J_{p,\d}(x)=\hat J_p(x),\qquad \forall\ x\in X.$$
}

\proof  (a) Follows by using Eq.\ \limdeltazero\ for $x\in \widehat X_p$, and by taking $w_{p,\d}(x)=0$ for $x\notin \widehat X_p$.
\smskip
\pn (b) By Prop.\ \propnegdp(e), there exists an $\e$-optimal policy $\p_\e$ for the $p$-$\d$-perturbed problem, i.e.,  $J_{\p_\e,p,\d}(x)\le \hat J_{p,\d}(x)+\e$ for all $x\in X$. Moreover, for  $x\in \widehat X_p$ we have $\hat J_{p,\d}(x)<\infty$, so $J_{\p_\e,p,\d}(x)<\infty$. Hence $\p_\e$ is $p$-stable from all $x\in \widehat X_p$, and we have $\hat J_p\le J_{\p_\e}$. 
Using also Eq.\ \approxform, we have  for all $\d>0$, $\e>0$, $x\in X$, and $\p\in \P_{p,x}$,
$$\hat J_p(x)-\e\le J_{\p_\e}(x)-\e\le J_{\p_\e,p,\d}(x)-\e\le  \hat J_{p,\d}(x)\le J_{\p,p,\d}(x)= J_{\p}(x)+w_{\p,p,\d}(x),$$
where $\lim_{\d\downarrow0}w_{\p,p,\d}(x)=0$ for all $x\in X$. 
By taking the limit as $\e\downarrow 0$, we obtain for all $\d>0$ and $\p\in \P_{p,x}$,
$$\hat J_p(x)\le \hat J_{p,\d}(x)\le J_{\p}(x)+w_{\p,p,\d}(x),\qquad \forall\ x\in X.$$
 By taking the limit as $\d\downarrow0$ and then the infimum over all  $\p\in \P_{p,x}$, we have
$$\hat J_p(x)\le \lim_{\d\downarrow0}\hat J_{p,\d}(x)\le \inf_{\p\in \P_{p,x}}J_\p(x)=\hat J_p(x),\qquad \forall\ x\in X,$$
from which the result follows.
\qed

We now consider approximately optimal policies. Given any $\e>0$, by Prop.\ \propnegdp(e), there exists an $\e$-optimal policy for the $p$-$\d$-perturbed problem, i.e., a policy $\p$ such that 
$J_\p(x)\le \hat J_{p,\d}(x)+\e$ for all $x\in X$. We address the question whether there exists a $p$-stable policy $\p$ that is $\e$-optimal for the restricted optimization over $p$-stable policies, i.e., a policy $\p$ that is $p$-stable simultaneously from all $x\in X_p$, (i.e., $\p\in \P_p$) and satisfies
$$J_\p(x)\le \hat J_p(x)+\e,\qquad\forall\ x\in X.$$
We refer to such a policy as a {\it $p$-$\e$-optimal policy\/}.

\xdef\propepsilonopt{\propn}\propnum\show{myproposition}

\texshopbox{\proposition{\propepsilonopt\ (Existence of $p$-$\e$-Optimal Policy):}Let $p$ be a given forcing function and $\d>0$. For every $\e>0$, a policy $\p$ that is $\e$-optimal for the $p$-$\d$-perturbed problem is $p$-$\e$-optimal, and hence belongs to $\P_p$. }

\proof For any $\e$-optimal policy $\p_\e$ for the $p$-$\d$-perturbed problem,  we have 
$$J_{\p_\e,p,\d}(x)\le \hat J_{p,\d}(x)+\e<\infty,\qquad \forall\ x\in \widehat X_p.$$
This implies that $\p_\e\in \P_p$. Moreover, for all sequences $\{x_k\}$ generated  
from initial state-policy pairs $(\p,x_0)$ with $x_0\in \widehat X_p$ and $\p\in\P_{p,x_0}$, we have
$$J_{\p_\e}(x_0)\le J_{\p_\e,p,\d}(x_0)\le \hat J_{p,\d}(x_0)+\e\le J_\p(x_0)+\d\sum_{k=0}^\infty p(x_k)+\e.$$
Taking the limit as $\d\downarrow0$ and using the fact $\sum_{k=0}^\infty p(x_k)<\infty$ (since $\p\in \P_{p,x_0}$), we obtain
$$J_{\p_\e}(x_0)\le J_\p(x_0)+\e,\qquad \forall\ x_0\in \widehat X_p,\ \p\in\P_{p,x_0}.$$
By taking infimum over  $\p\in\P_{p,x_0}$, it follows that 
$$J_{\p_\e}(x_0)\le \hat J_p(x_0)+\e,\qquad \forall\ x_0\in \widehat X_p,$$
which in view of the fact $J_{\p_\e}(x_0)= \hat J_p(x_0)=\infty$ for $x_0\notin \widehat X_p,$ implies that $\p_\e$ is $p$-$\e$-optimal. \qed

Note that the preceding proposition implies that
$$\hat J_p(x)=\inf_{\p\in \Pi_p}J_{\p}(x),\qquad\forall\ x\in X,\xdef\epsopt{\lab}\eqnum\show{oneo}$$
which is a stronger statement than the definition $\hat J_p(x)=\inf_{\p\in \Pi_{p,x}}J_{\p}(x)$ for all $x\in X$. However, it can be shown through examples that there may not exist a restricted-optimal $p$-stable policy, i.e., a $\p\in\P_p$ such that $J_\p=\hat J_p$, even if there exists an optimal policy for the original problem. One such example is the one-dimensional linear-quadratic problem of Example \examplelq\ for the case where $p=p^+$. Then, there exists a unique linear stable policy that attains the restricted optimal cost $J^+(x)$ for all $x$ (cf.\ Fig.\ \figlinquad), but this policy is not terminating. Note also that there may not exist a {\it stationary} $p$-$\e$-optimal policy, since generally in undiscounted nonnegative optimal control problems there may not exist a stationary $\e$-optimal policy, as is well-known (for an example, see [Ber13], following Prop.\ 4.3.2).

Our next proposition is preliminary for our main result. It involves the set of functions $S_p$ given by
$$\eqalign{S_p=\Big\{J\in{\cal J} \ \big |\ & 
J(x_k)\to0 \hbox{ for all sequences $\{x_k\}$ generated  
from}\cr
&\hbox{initial state-policy pairs $(\p,x_0)$ with $x_0\in X$ and $\p\in\P_{p,x_0}$} \Big\}.\cr}\xdef\widehatj{\lab}\eqnum\show{oneo}$$
In words, $S_p$ is the set of functions in ${\cal J}$ whose value is asymptotically driven to 0 by all the policies that are $p$-stable starting from some $x_0\in X$ (and thus have the character of Lyapounov functions for these policies).\old{\footnote{\dag}{\ninepoint Note here that the state-policy pairs $(\p,x_0)$ with $x_0\in X$ and $\p\in\P_{p,x_0}$ are the same as the ones with $x_0\in \widehat X_p$ and $\p\in\P_{p,x_0}$, since $\P_{p,x_0}=\emptyset$ for $x_0\notin \widehat X_p$.}}

Note that $S_p$ contains $\hat J_p$ and $\hat J_{p,\d}$  for all $\d>0$ [cf.\ Eq.\ \limprophatj]. Moreover, $S_p$  contains all functions $J$ such that $0\le J\le h(\hat J_{p,\d})$ 
for some $\d>0$ and function $h:X\mapsto X$ such that $h(J)\to0$ as $J\to0$. For example $S_p$  contains all  $J$ such that $0\le J\le c\,\hat J_{p,\d}$ for some $c>0$ and $\d>0$. 

We summarize the preceding discussion in the following proposition, which also shows uniqueness of solution (within $S_p$) of Bellman's equation for the $p$-$\d$-perturbed problem. The significance of this is that  the $p$-$\d$-perturbed problem, which can be solved more reliably than the original problem (including by VI methods), can yield a close approximation to $\hat J_p$ for sufficiently small $\d$ [cf.\ Prop.\ \propdeltaper(b)].

\xdef\propdeltaperbel{\propn}\propnum\show{myproposition}

\texshopbox{\proposition{\propdeltaperbel:} Let $p$ be a forcing function and $\d>0$. The function $\hat J_{p,\d}$ belongs to the set $S_p$, and is the unique solution within $S_p$ of Bellman's equation for the $p$-$\d$-perturbed problem,
$$\hat J_{p,\d}(x)=\inf_{u\in U(x)}\Big\{g(x,u)+\d p(x)+\hat J_{p,\d}\big(f(x,u)\big)\Big\},\qquad x\in X.\xdef\deltabeleq{\lab}\eqnum\show{oneo}$$
Moreover,  $S_p$ contains $\hat J_p$ and  all functions $J$ satisfying 
$$0\le J\le h(\hat J_{p,\d})$$ 
for some $h:X\mapsto X$ with $h(J)\to0$ as $J\to0$. 
}

\proof We have $\hat J_{p,\d}\in S_p$ and $\hat J_p\in S_p$ by Eq.\ \limprophatj, as noted earlier. We also have that $\hat J_{p,\d}$ is a solution of Bellman's equation \deltabeleq\ by Prop.\ \propnegdp(a).
To show that $\hat J_{p,\d}$ is the unique solution within $S_p$, let $\tl J\in S_p$ be another solution, so that using also Prop.\ \propnegdp(a), we have 
$$\hat J_{p,\d}(x)\le \tl J(x)\le g(x,u)+\d p(x)+\tl J \big(f(x,u)\big),\qquad \forall\ x\in X,\ u\in U(x).\xdef\deltatlineq{\lab}\eqnum\show{oneo}$$
Fix $\e>0$, and let $\p=\{\m_{0},\m_{1},\ldots\}$ be an $\e$-optimal policy for the $p$-$\d$-perturbed problem.
By repeatedly applying the preceding relation, we have for any $x_0\in \widehat X_p$,
$$\hat J_{p,\d}(x_0)\le \tl J(x_0)\le \tl J(x_k)+\d\sum_{m=0}^{k-1}p(x_m)+\sum_{m=0}^{k-1}g\big(x_m,\m_{m}(x_m)\big),\qquad \forall\ k\ge1,\xdef\relone{\lab}\eqnum\show{oneo}$$
where $\{x_k\}$ is the state sequence generated starting from $x_0$ and using $\p$. We have $\tl J(x_k)\to0$ (since $\tl J\in S_p$ and $\p\in\P_p$ by Prop.\ \propepsilonopt), so that
$$\lim_{k\to\infty}\lf\{\tl J(x_k)+\d\sum_{m=0}^{k-1}p(x_m)+\sum_{m=0}^{k-1}g\big(x_m,\m_{m}(x_m)\big)\ri\}=J_{\p,\d}(x_0)\le \hat J_{p,\d}(x_0)+\e.\xdef\reltwo{\lab}\eqnum\show{oneo}$$
By combining Eqs.\ \relone\ and \reltwo, we obtain
$$\hat J_{p,\d}(x_0)\le \tl J(x_0)\le \hat J_{p,\d}(x_0)+\e,\qquad \forall\ x_0\in \widehat X_p.$$
By letting $\e\to0$, it follows that $\hat J_{p,\d}(x_0)= \tl J(x_0)$ for all $x_0\in \widehat X_p$. Also for $x_0\notin \widehat X_p$, we have 
$\hat J_{p,\d}(x_0)= \tl J(x_0)=\infty$
 [since $\hat J_{p,\d}(x_0)=\infty$ for $x_0\notin \widehat X_p$ and $\hat J_{p,\d}\le \tl J$, cf.\ Eq.\ \deltatlineq]. Thus $\hat J_{p,\d}= \tl J$, proving that $\hat J_{p,\d}$ is the unique solution of the Bellman Eq.\ \deltabeleq\ within $S_p$. 
\qed

We next show that $\hat J_p$ is the unique solution of Bellman's equation within the set of functions
$${\cal W}_p=\{J\in S_p\mid \hat J_p\le J\},\xdef
\setcalw{\lab}\eqnum\show{oneo}$$
and that the VI algorithm yields $\hat J_p$ in the limit for any initial $J_0\in {\cal W}_p$.

\xdef\proptaffmonft{\propn}\propnum\show{myproposition}

\texshopbox{\proposition{\proptaffmonft:} Let $p$ be a given forcing function. Then:
\nitem{(a)} $\hat J_p$ is the unique solution of Bellman's equation
$$J(x)=\inf_{u\in U(x)}\Big\{g(x,u)+J\big(f(x,u)\big)\Big\},\qquad x\in X,\xdef\bellmaneq{\lab}\eqnum\show{oneo}$$
 within the set ${\cal W}_p$ of Eq.\ \setcalw.
\nitem{(b)} ({\it VI Convergence\/})  If $\{J_k\}$ is the sequence generated by the VI algorithm \vieq\ starting with some 
 $J_0\in {\cal W}_p$, then $J_k\to\hat J_p$.
\nitem{(c)}  ({\it Optimality Condition\/})  If $\hat \m$ is a $p$-stable stationary policy and 
$$\hat \m(x)\in \argmin_{u\in U(x)}\big\{g(x,u)+\hat J_p\big(f(x,u)\big)\big\},\qquad \forall\ x\in X ,\xdef
\optcond{\lab}\eqnum\show{oneo}$$
then $\hat \m$ is optimal over the set of $p$-stable policies. Conversely, if $\hat \m$ is optimal within the set of $p$-stable policies, then it satisfies the preceding condition \optcond.
}

\proof (a), (b) We first show that $\hat J_p$ is a solution of Bellman's equation. Since $\hat J_{p,\d}$ is a solution [cf.\ Prop.\ \propdeltaperbel] and $\hat J_{p,\d}\ge \hat J_p$ [cf.\ Prop.\ \propdeltaper(b)], we have for all $\d>0$,
$$\eqalignno{\hat J_{p,\d}(x)&=\inf_{u\in U(x)}\Big\{g(x,u)+\d p(x)+\hat J_{p,\d}\big(f(x,u)\big)\Big\}\cr
&\ge \inf_{u\in U(x)}\Big\{g(x,u)+\hat J_{p,\d}\big(f(x,u)\big)\Big\}\cr
&\ge \inf_{u\in U(x)}\Big\{g(x,u)+\hat J_p\big(f(x,u)\big)\Big\}.\cr}$$
By taking the limit as $\d\downarrow0$ and using the fact $\lim_{\d\downarrow0}\hat J_{p,\d}=\hat J_p$ [cf.\ Prop.\ \propdeltaper(b)], we obtain
$$\hat J_p(x)\ge \inf_{u\in U(x)}\Big\{g(x,u)+\hat J_p\big(f(x,u)\big)\Big\},\qquad \forall\ x\in X.\xdef
\forineq{\lab}\eqnum\show{oneo}$$

For the reverse inequality, let $\{\d_m\}$ be a sequence with $\d_m\downarrow0$.  From Prop.\ \propdeltaperbel, we have for all $m$, $x\in X$, and $u\in U(x)$,
$$g(x,u)+\d_m p(x)+\hat J_{p,\d_m}\big(f(x,u)\big)
\ge \inf_{v\in U(x)}\Big\{g(x,v)+\d_m p(x)+\hat J_{p,\d_m}\big(f(x,v)\big)\Big\}=\hat J_{p,\d_m}(x).$$
Taking the limit as $m\to\infty$, and using the fact $\lim_{\d_m\downarrow0}\hat J_{p,\d_m}=\hat J_p$ [cf.\ Prop.\ \propdeltaper(b)], we have
$$g(x,u)+\hat J_p\big(f(x,u)\big)\ge \hat J_p(x),\qquad \forall\ x\in X,\ u\in U(x),$$
so that 
$$\inf_{u\in U(x)}\Big\{g(x,u)+\hat J_p\big(f(x,u)\big)\Big\}\ge \hat J_p(x),\qquad \forall\ x\in X.\xdef
\backineq{\lab}\eqnum\show{oneo}$$
By combining Eqs.\ \forineq\ and \backineq, we see that $\hat J_p$ is a solution of Bellman's equation. We also have $\hat J_p\in S_p$ by Prop.\ \propdeltaperbel, implying that $\hat J_p\in {\cal W}_p$ and proving part (a) except for the uniqueness assertion.

We will now prove part (b). Let $\p=\{\m_0,\m_1,\ldots\}\in\P_p$ [which is nonempty by Prop.\ \propdeltaper(c)], and for $x_0\in \widehat X_p$, let $\{x_k\}$ be the generated sequence starting from $x_0$ and using $\p$. We have $J_0(x_k)\to0$  since $J_0\in S_p$. Since from the definition of the VI sequence $\{J_k\}$, we have
$$J_k(x)\le g(x,u)+J_{k-1}\big(f(x,u)\big),\qquad \forall\ x\in X,\ u\in U(x),\ k=1,2,\ldots,$$
it follows that
$$J_k(x_0)\le J_0(x_k)+\sum_{m=0}^{k-1}g\big(x_m,\m_m(x_m)\big).$$
By taking limit as $k\to\infty$ and using the fact $J_0(x_k)\to0$, it follows that $\limsup_{k\to\infty}J_k(x_0)\le J_\p(x_0)$. By taking the infimum over all $\p\in\P_p$, we obtain 
$\limsup_{k\to\infty}J_k(x_0)\le \hat J_p(x_0).$ Conversely, since $\hat J_p\le J_0$ and $\hat J_p$ is a solution of Bellman's equation (as shown earlier), it follows by induction that $\hat J_p\le J_k$ for all $k$. Thus 
$\hat J_p(x_0)\le \liminf_{k\to\infty}J_k(x_0),$
 implying that $J_k(x_0)\to\hat J_p(x_0)$ for all $x_0\in \widehat X_p$. We also have $\hat J_p\le J_k$ for all $k$, so that $\hat J_p(x_0)=J_k(x_0)=\infty$ for all $x_0\notin \widehat X_p$. This completes the proof of part (b).
Finally, since $\hat J_p\in {\cal W}_p$ and $\hat J_p$ is a solution of Bellman's equation, part (b) implies the uniqueness assertion of part (a).

\smskip
\pn (c) If $\m$ is $p$-stable and Eq.\ 
\optcond\ holds, then 
$$\hat J_p(x)=g\big(x,\m(x)\big)+\hat J_p\big(f(x,\m(x))\big),\qquad x\in X.$$
By Prop.\ \propnegdp(b), this implies that  
$J_\m\le \hat J_p$, so 
$\m$ is optimal over the set of $p$-stable policies. Conversely, assume that $\m$ is $p$-stable and $J_\m=\hat J_p$. Then by Prop.\ \propnegdp(b), we have
$$\hat J_p(x)=g\big(x,\m(x)\big)+\hat J_p\big(f(x,\m(x))\big),\qquad x\in X,$$
and since [by part (a)] $\hat J_p$ is a solution of Bellman's equation,
$$\hat J_p(x)=\inf_{u\in U(x)}\big\{g(x,u)+\hat J_p\big(f(x,u)\big)\big\},\qquad x\in X.$$
Combining the last two relations, we obtain Eq.\ \optcond.  \qed

We now consider the special case where $p$ is equal to the function $p^+(x)=1$ for all $x\ne t$ [cf.\ Eq.\ \terminatingp]. The set of $p^+$-stable policies from $x$ is $\P^+_x$, the set of terminating policies from $x$, and  $J^+(x)$ is the corresponding restricted optimal cost, 
$$J^+(x)=\hat J_{p^+}(x)=\inf_{\p\in\P_x^+}J_\p(x)=\inf_{\p\in\P^+}J_\p(x),\qquad x\in X,$$
[the last equality follows from Eq.\ \epsopt]. In this case, the set $S_{p^+}$ of Eq.\ \widehatj\ is the entire set ${\cal J}$, since for all $J\in{\cal J}$ and all sequences $\{x_k\}$ generated  
from initial state-policy pairs $(\p,x_0)$ with $x_0\in X$ and $\p$ terminating from $x_0$, we have $J(x_k)=0$ for $k$ sufficiently large. Thus, the set ${\cal W}_{p^+}$ of Eq.\ 
\setcalw\ is 
$${\cal W}_{p^+}=\{J\in {\cal J}\mid J^+\le J\}.\xdef
\wplus{\lab}\eqnum\show{oneo}$$
By specializing to the case $p=p^+$ the result of Prop.\ \proptaffmonft, we obtain the following proposition, which makes a stronger assertion than Prop.\ \proptaffmonft(a), namely that $J^+$ is the largest solution of Bellman's equation within ${\cal J}$ (rather than the smallest solution within ${\cal W}_{p^+}$).

\xdef\propbeleqterm{\propn}\propnum\show{myproposition}

\texshopbox{\proposition{\propbeleqterm:} 
\nitem{(a)} $J^+$ is the largest solution of the Bellman equation \bellmaneq\ within ${\cal J}$,
i.e., if $\tl J\in {\cal J}$ is a solution of Bellman's equation, then $\tl J\le J^+$.
\nitem{(b)} ({\it VI Convergence\/})  If $\{J_k\}$ is the sequence generated by the VI algorithm \vieq\ starting with some 
 $J_0\in {\cal J}$ with $J_0\ge J^+$, then $J_k\to J^+$.
\nitem{(c)}  ({\it Optimality Condition\/})  If $\m^+$ is a terminating stationary policy and 
$$\m^+(x)\in \argmin_{u\in U(x)}\big\{g(x,u)+J^+\big(f(x,u)\big)\big\},\qquad \forall\ x\in X ,\xdef\optcondterm{\lab}\eqnum\show{oneo}$$
then $\m^+$ is optimal over the set of terminating policies. Conversely, if $\m^+$ is optimal within the set of terminating policies, then it satisfies the preceding condition \optcondterm.
}

\proof In view of Prop.\ \proptaffmonft\ and the expression \wplus\ for ${\cal W}_{p^+}$, we only need to show that $\tl J\le J^+$ for every solution $\tl J\in {\cal J}$ of Bellman's equation. Indeed, let $\tl J$ be such a solution. We have $\tl J(x_0)\le J^+(x_0)$ for all $x_0$ with $J^+(x_0)=\infty$, so in order to show that $\tl J\le J^+$, it will suffice to show that for every $(\p,x_0)$ with $\p\in \P^+_{x_0}$, we have $\tl J(x_0)\le J_\p(x_0)$. 
Indeed, consider $(\p,x_0)$ with $\p\in \P^+_{x_0}$, and let $\{x_0,\ldots,x_k,t\}$ be the terminating state sequence generated starting from $x_0$ and using $\p$. Since $\tl J$ solves Bellman's equation, we have
$$\eqalign{\tl J(x_m)&\le g\big(x_m,\m_m(x_m)\big)+\tl J(x_{m+1}),\qquad m=0,\ldots,k-1,\cr
\tl J(x_k)&\le g\big(x_k,\m_k(x_k)\big).\cr}$$
By adding these relations, we obtain 
$$\tl J(x_0)\le \sum_{m=0}^k g\big(x_m,\m_m(x_m)\big)=J_\p(x_0),\qquad \forall\ (\p,x_0)\hbox{ with } \p\in \P^+_{x_0},$$
and by taking the infimum of the right side over $\p\in \P^+_{x_0}$, we obtain $\tl J(x_0)\le J^+(x_0)$. \qed

\xdef \figbellmansol{\figr}\figrnum\show{myfigure}

{
\topinsert
\centerline{\hskip0pc\includegraphics[width=4.5in]{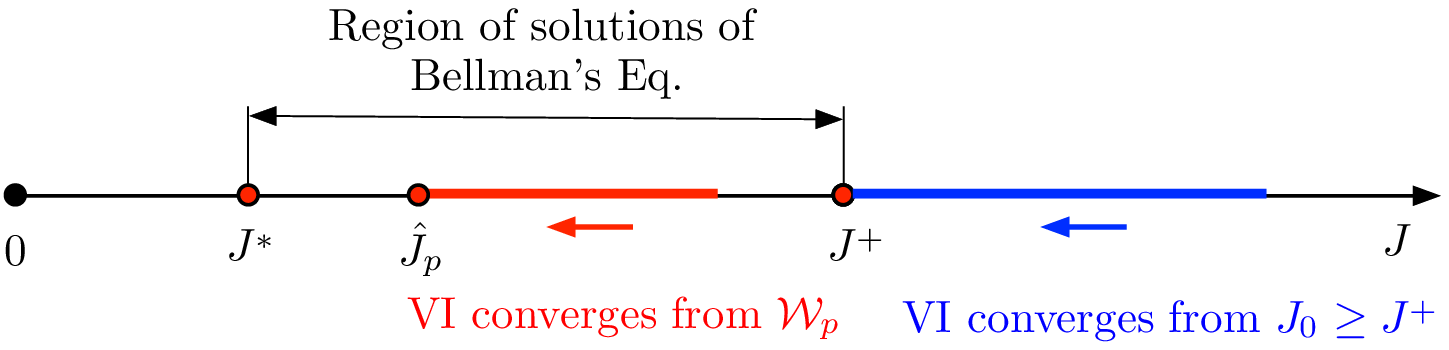}}
\fig{-0.2pc}{\figbellmansol} {Illustration of the  solutions of Bellman's equation. The smallest and the largest solutions are $J^*$ and $J^+$, respectively. The VI algorithm converges to $J^+$ starting from any $J_0\in{\cal J}$ with $J_0\ge J^+$, and it converges to $\skew6\hat J_p$ starting from any $J_0\in {\cal W}_p$.}
\endinsert
}

\xdef \figlyapounov{\figr}\figrnum\show{myfigure}

 We illustrate Props.\ \proptaffmonft\ and \propbeleqterm\ in Figs.\ \figbellmansol\ and \figlyapounov. In particular, each forcing function $p$ delineates the set of initial functions ${\cal W}_p$ from which VI converges to $\hat J_p$. The function $\hat J_p$ is the minimal element of ${\cal W}_p$. Moreover, we have 
${\cal W}_p\cap {\cal W}_{p'}=\emptyset\quad \hbox {if}\quad \hat J_p\ne \hat J_{p'},$ in view of the VI convergence result of Prop.\ \proptaffmonft(b).

{
\topinsert
\centerline{\hskip0pc\includegraphics[width=4.5in]{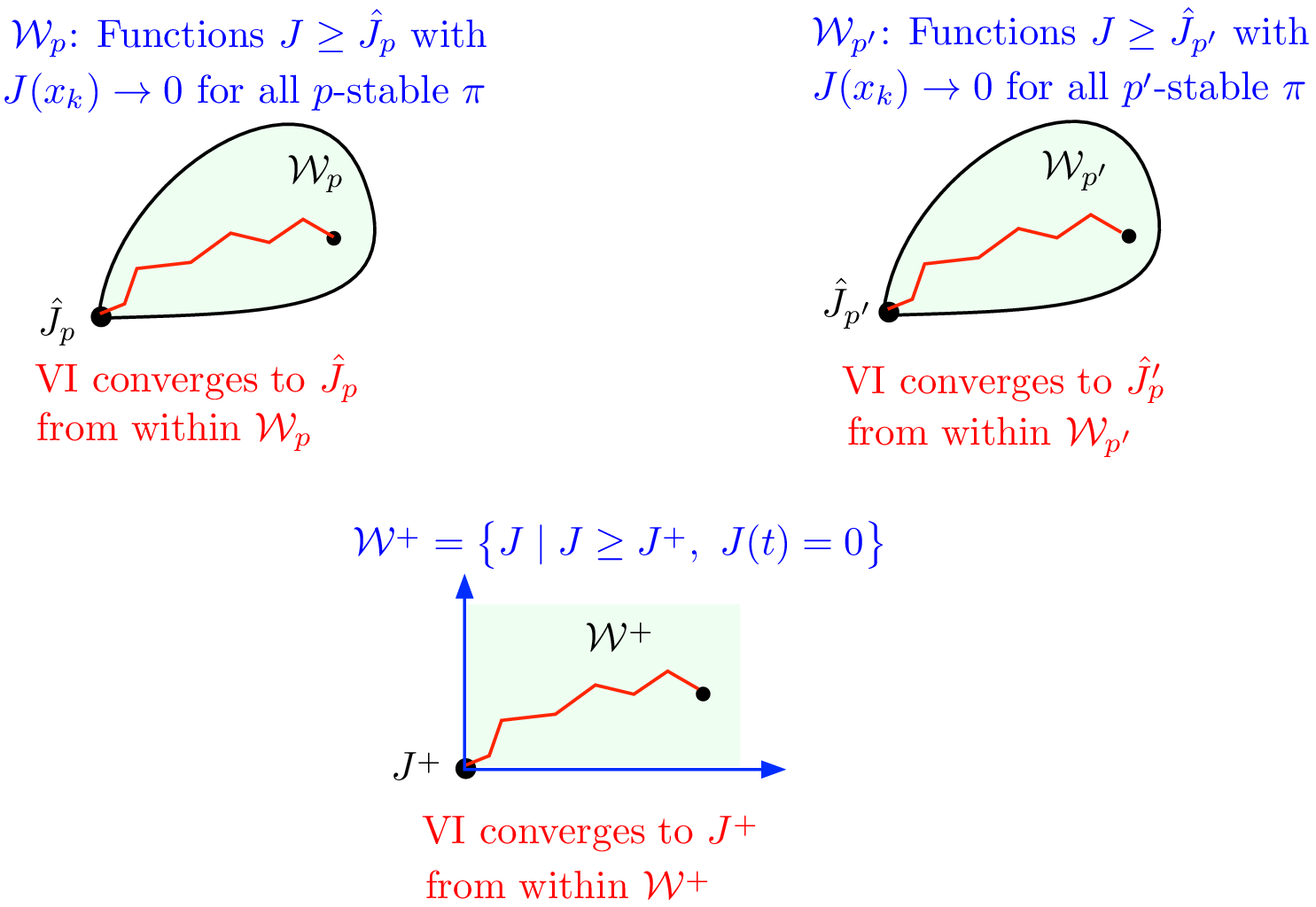}}
\fig{-0.2pc}{\figlyapounov} {Illustration of the VI convergence results of Prop.\ \proptaffmonft\ and \propbeleqterm. Each $p$ defines the set of initial functions ${\cal W}_p$ from which VI converges to $\skew6\hat J_p$ from above. For two forcing functions $p$ and $p'$, we have ${\cal W}_p\cap {\cal W}_{p'}=\emptyset$ if $\skew6\hat J_p\ne \skew6\hat J_{p'}$. 
}
\endinsert
}

Note that Prop.\ \propbeleqterm(b) implies that VI converges to $J^+$ starting from the particular initial condition 
$$J_0(x)=\cases{0&if $x=t$,\cr
\infty&if $x\ne t$.\cr}\eqnum\show{oneo}$$
For this choice of $J_0$, the value $J_k(x)$ generated by VI is the optimal cost that can be achieved starting from $x$ subject to the constraint that $t$ is reached in $k$ steps or less. 

Suppose now that the set of terminating policies is sufficient in the sense that it can achieve the same optimal cost as the set of all policies, i.e., $J^+=\jstar$.  Then, from Prop.\ \propbeleqterm, it follows that $\jstar$ is the unique solution of Bellman's equation within ${\cal J}$, and the VI algorithm converges to $\jstar$ from above, i.e., starting from any $J_0\in {\cal J}$ with $J_0\ge \jstar$. Under additional conditions, such as finiteness of $U(x)$ for all $x\in X$ [cf.\ Prop.\ \propnegdp(d)], VI converges to $\jstar$ starting from any $J_0\in{\cal J}$. 

Examples of problems where terminating policies are sufficient include linear-quadratic problems under the classical conditions of controllability and observability (cf.\ Example \examplelq), and finite-node deterministic shortest path problems with all cycles having positive length. Note that in the former case, despite the fact $J^+=\jstar$, there is no optimal terminating policy, since the only optimal policy is a linear policy that drives the system to the origin asymptotically, but not in finite time.

Let us illustrate the results of this section with two examples.

\xdef\examplelqminenergy{\exampl}\examplnum\show{examplo}

 \beginexample{\examplelqminenergy\ (Minimum Energy Stable Control of \hfill\break Linear Systems)}Consider the linear-quadratic problem of Example \examplelq. We assume that the pair $(A,B)$ is stabilizable. However, we are making no assumptions on the state weighting matrix $Q$ other than positive semidefiniteness, so the detectability assumption may not be satisfied. This includes the case $Q=0$, when $J^*(x)\equiv0$. In this case an optimal policy is $\m^*(x)\equiv0$, which may not be stable, yet the problem of finding a stable policy that minimizes the ``control  energy" (a cost that is quadratic on the control with no penalty on the state) among all stable policies  is meaningful.
 
 We consider the forcing function 
 $$p(x)=\|x\|^2,$$
  so the $p$-$\d$-perturbed problem satisfies the detectability condition and from classical results, $\skew4\hat J_{p,\d}$ is a positive definite quadratic function $x'P_\d x$, where $P_\d$ is the unique solution of the $\d$-perturbed Riccati equation
 $$P_\d = A'\bl( P_\d-P_\d B(B'P_\d B + R)^{-1}B'P_\d\br) A+Q+\d I,\xdef
\pertriccati{\lab}\eqnum\show{oneo}$$
within the class of positive semidefinite matrices. By Prop.\ \propdeltaper, we have $\skew4\hat J_p(x)=x'\skew4 \hat Px$, where $\skew4\hat P=\lim_{\d\downarrow0}P_\d$ is positive semidefinite, and solves the (unperturbed) Riccati equation
$$P= A'\bl( P-P B(B'P B + R)^{-1}B'P\br) A+Q.$$ 
Moreover, by Prop.\ \proptaffmonft(a), $\skew4\hat P$ is the largest solution among positive semidefinite matrices, since all positive semidefinite quadratic functions belong to the set $S_p$ of Eq.\ \widehatj.  By Prop.\ \proptaffmonft(c), any stable stationary policy $\skew3\hat \m$ that is optimal among the set of stable policies must satisfy the optimality condition
 $$\skew3\hat \m(x)\in \argmin_{u\in \re^m }\big\{u'Ru+(Ax+Bu)'\skew4\hat P(Ax+Bu)\big\},\qquad \forall\ x\in \rn,$$
 [cf.\ Eq.\ \optcond],  or equivalently, by setting the gradient of the minimized expression to 0,
 $$(R+B'\skew4\hat P B)\skew3\hat\m(x)=-B'\skew4\hat PAx.\xdef
\optcondric{\lab}\eqnum\show{oneo}$$ 
We may solve Eq.\ \optcondric, and check if any of its solutions $\skew3\hat \m$ is $p$-stable; if this is so, $\skew3\hat \m$ is optimal within the class of $p$-stable policies. Note, however, that in the absence of additional conditions, it is possible that some policies $\skew3\hat \m$ that solve Eq.\ \optcondric\ are $p$-unstable. 

In the case where the pair $(A,B)$ is not stabilizable, the $p$-$\d$-perturbed cost function $\skew6\hat J_{p,\d}$ need not be real-valued, and the $\d$-perturbed Riccati equation \pertriccati\ may not have any solution (consider for example the case where $n=1$, $A=2$, $B=0$, and $Q=R=1$). Then, Prop.\ \propbeleqterm\ still applies, but the preceding analytical approach needs to be modified.
\endexample

As noted earlier, the Bellman equation may have multiple solutions corresponding to different forcing functions $p$, with each solution being unique within the corresponding set ${\cal W}_p$ of Eq.\ \setcalw, consistently with Prop.\ \proptaffmonft(a). The following is an illustrative example.

\xdef \figstopping{\figr}\figrnum\show{myfigure}

\xdef \figstopscalar{\figr}\figrnum\show{myfigure}

\xdef\exampleperturb{\exampl}\examplnum\show{examplo}

 \beginexample{\exampleperturb\ (An Optimal Stopping Problem)}Consider an optimal stopping problem where the state space $X$ is $\rn$. We identify the destination with the origin of $\rn$, i.e., $t=0$. At each $x\ne 0$, we may either stop (move to the origin) at a cost $c>0$, or move to state $\g x$ at cost $\|x\|$, where $\g$ is a scalar with $0<\g<1$; see Fig.\ \figstopping.\footnote{\dag}{\ninepoint In this example, the salient feature of the policy that never stops is that it drives the system asymptotically to the destination according to an equation of the form $x_{k+1}=f(x_k)$, where $f$ is a contraction mapping. The example admits generalization to the broader class of optimal stopping problems where the policy that never stops has this property. For simplicity in illustrating our main point, we consider here the special case where $f(x)=\g x$ with $\g\in(0,1)$.} Thus the Bellman equation  has the form
 $$J(x)=\cases{\min\big\{c,\,\|x\|+J(\g x)\big\}&if $x\ne0$,\cr
 0&if $x=0$.\cr}\eqnum\show{threeo}$$

{
\topinsert
\centerline{\hskip0pc\includegraphics[width=2.2in]{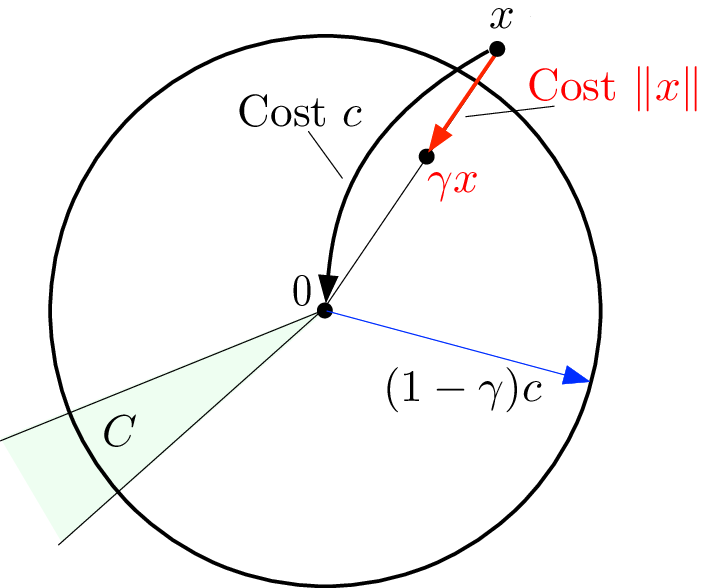}}
\fig{-0.2pc}{\figstopping} {Illustration of the stopping problem of Example \exampleperturb. The optimal policy is to stop outside the sphere of radius $(1-\g)c$ and to continue otherwise. Each cone $C$ of the state space defines a different solution $\skew6\hat J_p$ of Bellman's equation, with $\skew6\hat J_p(x)=c$ for all nonzero $x\in C$, and a corresponding region of convergence of the VI algorithm.}
\endinsert
}

Let us consider first the forcing function
 $$p(x)=\|x\|.$$
Then it can be verified that all policies are $p$-stable. We have
 $$J^*(x)=\skew6\hat J_p(x)=\min\lf\{c,\,{1\over 1-\g}\|x\|\ri\},\qquad\forall\ x\in\rn,$$
and the optimal cost function of the corresponding $p$-$\d$-perturbed problem is
 $$\skew6\hat J_{p,\d}(x)=\min\lf\{c+\d\|x\|,\,{1+\d\over 1-\g}\|x\|\ri\},\qquad\forall\ x\in\rn.$$
 Here the set $S_p$ of Eq.\ \widehatj\ is given by
 $$S_p=\Big\{J\in {\cal J}\mid \lim_{x\to0}J(x)=0\Big\},$$
and the corresponding set ${\cal W}_p$ of Eq.\ \setcalw\ is given by
$${\cal W}_p=\Big\{J\in {\cal J}\mid J^*\le J,\ \lim_{x\to0}J(x)=0\Big\}.$$

Let us consider next the forcing function
 $$p^+(x)=\cases{1&if $x\ne0$,\cr
 0&if $x=0$.\cr}$$
Then the $p^+$-stable policies are the terminating policies. Since stopping at some time and incurring the cost $c$ is a requirement for a terminating policy, it follows that the optimal $p^+$-stable policy is to stop as soon as possible, i.e., stop at every state. The corresponding restricted optimal cost function is
 $$ J^+(x)=\cases{c&if $x\ne0$,\cr
 0&if $x=0$.\cr}$$
The optimal cost function of the corresponding $p^+$-$\d$-perturbed problem is
 $$\skew6\hat J_{p^+,\d}(x)=\cases{c+\d&if $x\ne0$,\cr
 0&if $x=0$,\cr}$$ 
 since in the $p^+$-$\d$-perturbed problem it is again optimal to stop as soon as possible, at cost $c+\d$.  Here the set $S_{p^+}$ is equal to ${\cal J},$ and the corresponding set ${\cal W}_{p^+}$ is equal to 
 $\big\{J\in {\cal J}\mid J^+\le J\big\}.$
 
However, there are infinitely many additional solutions of Bellman's equation between the largest and smallest solutions $J^*$ and $J^+$. For example, when $n>1$, functions $J\in{\cal J}$ such that $J(x)=J^*(x)$ for $x$ in some cone and $J(x)=J^+(x)$ for $x$ in the complementary cone are solutions; see Fig.\ \figstopping.  There is also a corresponding infinite number of regions of convergence ${\cal W}_{p}$ of VI. Also VI  converges to $J^*$ starting from any $J_0$ with $0\le J_0\le J^*$ [cf.\ Prop.\ \propnegdp(d)]. Figure \figstopscalar\ illustrates additional  solutions of Bellman's equation of a different character.
\endexample

\topinsert
\centerline{\hskip0pc\includegraphics[width=5.5in]{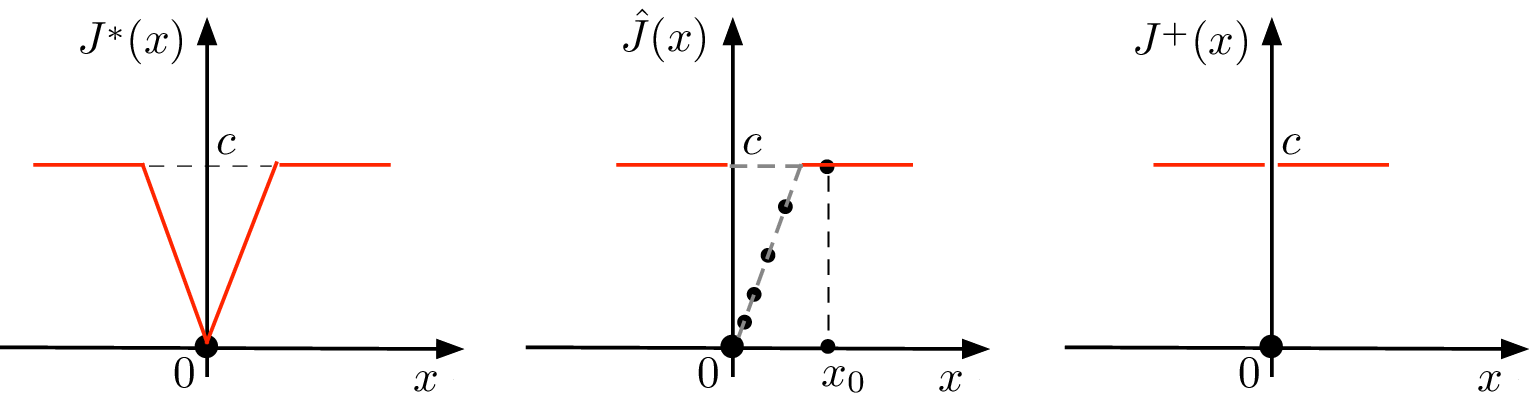}}
\fig{-0.0pc}{\figstopscalar} {Illustration of three solutions of Bellman's equation in the one-dimensional case ($n=1$) of the stopping problem of Example \exampleperturb. 
The solution in the middle is specified by a scalar $x_0>0$, and has the form
$$\skew6\hat J(x)=\cases{0&if $x=0$,\cr
{1\over 1-\g}|x|&if $0<x<(1-\g)c$ and $x=\g^k x_0$ for some $k\ge0$,\cr
c&otherwise.\cr}$$
}
\endinsert

\section{Policy Iteration Methods}

\pn Generally, the standard PI algorithm  \poleval, \polimprove\ produces unclear results under our assumptions. As an illustration, in the stopping problem of Example \exampleperturb, if PI is started with the policy that stops at every state, it repeats that policy, and this policy is not optimal  within the class of $p$-stable policies with respect to the forcing function $p(x)=\|x\|$. The following example provides an instance where the PI algorithm may converge to either an optimal or a strictly suboptimal policy.

\xdef\examplecounterpi{\exampl}\examplnum\show{examplo}
 
 \beginexample{\examplecounterpi\ (Counterexample for PI)}Consider the case $X=\{0,1\}$, $U(0)=U(1)=\{0,1\}$, and the destination is $t=0$. Let also
 $$f(x,u)=\cases{0&if $u=0$,\cr
 x&if $u=1$,\cr}\qquad g(x,u)=\cases{1&if $u=0,\ x=1$,\cr
 0&if $u=1\hbox{ or }x=0$.\cr}$$
 This is a shortest path problem where the control $u=0$ moves the state from $x=1$ to $x=0$ (the destination) at cost 1, while the control $u=1$ keeps the state unchanged at cost 0. The  policy $\m^*$ that keeps the state unchanged is the only optimal policy, with $J_{\m^*}(x)=J^*(x)=0$ for both states $x$. However, under any forcing function $p$ with $p(1)>0$, the policy $\skew3\hat \m$, which moves from state 1 to 0, is the only $p$-stable policy, and we have $J_{\skew3\hat \m}(1)=\skew5\hat J_p(1)=1$. The standard PI algorithm  \poleval, \polimprove\ if started with $\m^*$ will repeat $\m^*$. If this algorithm is started with $\skew4\hat \m$, it may generate $\m^*$ or it may repeat $\skew4\hat \m$, depending on how the policy improvement iteration is implemented. The reason is that for both $x$ we have
$$\skew4\hat \m(x)\in\arg\min_{u\in\{0,1\}}\Big\{g(x,u)+\skew5\hat J_p\big(f(x,u)\big)\Big\},$$
as can be verified with a straightforward calculation. Thus a rule for breaking a tie in the policy improvement operation is needed, but such a rule may not be obvious in general.
\endexample


Motivated by the preceding example, we consider several types of PI method that bypass the difficulty above either through assumptions or through modifications. We first consider a case where the PI algorithm is reliable. This is the case where the terminating policies are sufficient, in the sense that $J^+=J^*$.  

\subsection{Policy Iteration for the Case $\jstar=J^+$}

\pn The PI algorithm starts with a stationary policy $\m^0$, and generates a sequence of stationary policies $\{\m^k\}$ via a sequence of policy evaluations to obtain $J_{\m^k}$ from the equation
$$J_{\m^k}(x)=g\big(x,\m^k(x)\big)+J_{\m^k}\big(f\big(x,\m^k(x)\big)\big),\qquad x\in X,\xdef\poleval{\lab}\eqnum\show{oneo}$$
interleaved with policy improvements to obtain $\m^{k+1}$ from $J_{\m^k}$ according to
$$\m^{k+1}(x)\in \argmin_{u\in U(x)}\big\{g(x,u)+J_{\m^k}\big(f(x,u)\big)\big\},\qquad x\in X.\xdef
\polimprove{\lab}\eqnum\show{oneo}$$
We implicitly assume here that  the minimum in Eq.\ \polimprove\ is attained for each $x\in X$, which is true under some compactness condition on either $U(x)$ or the level sets of the function $g(x,\cdot)+J_k\big(f(x,\cdot)\big)$, or both.

\xdef\propdetpolicyit{\propn}\propnum\show{myproposition}

\texshopbox{\proposition{\propdetpolicyit: (Convergence of PI)} Assume that $\jstar=J^+$. Then the sequence $\{J_{\m^k}\}$ generated by the PI algorithm \poleval, \polimprove, satisfies $J_{\m^k}(x)\downarrow \jstar(x)$ for all $x\in X$. }

\proof
If $\m$ is a stationary policy and $\bar \m$ satisfies the policy improvement equation 
$$\bar \m(x)\in \argmin_{u\in U(x)}\big\{g(x,u)+J_{\m}\big(f(x,u)\big)\big\},\qquad x\in X,$$
[cf.\ Eq.\ \polimprove],  we have for all $x\in X$,
$$J_\mu(x)=g\big(x,\m(x)\big)+J_{\m}\big(f\big(x,\m(x)\big)\big)\ge \min_{u\in U(x)}\big\{g(x,u)+J_{\m}\big(f(x,u)\big)\big\}=g\big(x,\bar \m(x)\big)+J_{\m}\big(f\big(x,\bar \m(x)\big)\big),\xdef\imprineq{\lab}\eqnum\show{oneo}$$ 
where the first equality follows from by Prop.\ \propnegdp(b), and the second equality follows from the definition of $\bar \m$. Repeatedly applying this relation, we see that the sequence $\big\{\tl J_{k}(x_0)\big\}$ defined by 
$$\tl J_{k}(x_0)=J_\m(x_k)+\sum_{m=0}^{k-1}g\big(x_m,\bar \m(x_m)\big),\qquad k=1,2,\ldots,$$
is monotonically nonincreasing, where $\{x_k\}$ is the sequence generated starting from $x_0$ and using $\m$.
Moreover, from Eq.\ \imprineq\ we have 
$$J_\m(x_0)\ge \min_{u\in U(x_0)}\big\{g(x,u)+J_\m\big(f(x,u)\big)\big\}=\tl J_{1}(x_0)\ge\tl J_{k}(x_0),$$ 
for all $k$. This implies that
$$J_\m(x_0)\ge \min_{u\in U(x_0)}\big\{g(x,u)+J_\m\big(f(x,u)\big)\big\}\ge \lim_{k\to\infty}\tl J_k(x_0)\ge \lim_{k\to\infty}\sum_{m=0}^{k-1}g\big(x_m,\bar \m(x_m)\big)=J_{\bar \m}(x_0),$$
where the last inequality follows since $J_\m\ge0$. In conclusion, we have
$$J_\m(x)\ge \inf_{u\in U(x)}\big\{g(x,u)+J_{\m}\big(f(x,u)\big)\big\}\ge J_{\bar \m}(x),\qquad x\in X.\xdef\imprineqo{\lab}\eqnum\show{oneo}$$
Using $\m^k$ and $\m^{k+1}$ in place of $\m$ and $\bar \m$, we see that the sequence $\{J_{\m^k}\}$ generated by PI converges monotonically to some function $J_\infty\in E^+(X)$, i.e., $J_{\m^k}\downarrow J_\infty$. Moreover, from Eq.\ \imprineqo\ we have
$$J_\infty(x)\ge \inf_{u\in U(x)}\big\{g(x,u)+J_\infty\big(f(x,u)\big)\big\},\qquad x\in X,$$
as well as
$$g(x,u)+J_{\m_k}\big(f(x,u)\big)\ge J_\infty(x),\qquad x\in X,\ u\in U(x).$$
We now take the limit in the second relation as $k\to\infty$, then the infimum over $u\in U(x)$, and then combine with the first relation, to obtain
$$J_\infty(x)= \inf_{u\in U(x)}\big\{g(x,u)+J_\infty\big(f(x,u)\big)\big\},\qquad x\in X.$$
Thus $J_\infty$ is a solution of Bellman's equation, satisfying $J_\infty\ge \jstar$ (since $J_{\m^k}\ge \jstar$ for all $k$) and $J_\infty\in{\cal J}$ (since $J_{\m^k}\in{\cal J}$), so  by Prop.\ \propbeleqterm(a), it must satisfy $J_\infty=\jstar$. \qed

\vskip-1pc

\subsection{A Perturbed Version of Policy Iteration}

\pn 
We now consider a PI algorithm that does not require the condition $\jstar=J^+$. We will provide a  version of the PI algorithm that uses the forcing function $p$ and generates a sequence $\{\m^k\}$ of $p$-stable policies such that $J_{\m^k}\to \hat J_p$. In this section, the forcing function $p$ is kept fixed, and to simplify notation, we abbreviate $J_{\m,p,\d}$ with $J_{\m,\d}$. The following assumption requires that the algorithm generates $p$-stable policies exclusively, which can be quite restrictive. For example it is not satisfied for the problem of Example \examplecounterpi.

\xdef\assumptionpertpi{\assumptionn}\assumptionnum\show{myproposition}

\texshopbox{\assumption{\assumptionpertpi:} For each $\d>0$ there exists at least one $p$-stable stationary policy $\m$ such that $J_{\m,\d}\in S_p$. Moreover, given a $p$-stable stationary policy $\m$ and a scalar $\d>0$, every stationary policy $\ol \m$ such that 
$${\ol \m}(x)\in\arg\min_{u\in U(x)}\Big\{g(x,u)+J_{\m,\d}\big(f(x,u)\big)\Big\},\qquad \forall\ x\in X,\xdef
\polimprove{\lab}\eqnum\show{oneo}$$
is $p$-stable, and at least one such policy exists.
}

The perturbed version of the PI algorithm is defined as follows. 
Let $\{\d_k\}$ be a positive sequence with $\d_k\downarrow 0$, and let $\m^0$ be a $p$-stable policy that satisfies $J_{\m^0,\d_0}\in S_p$. One possibility is that $\m^0$ is an optimal policy for the $\d_0$-perturbed problem (cf.\ the discussion preceding Prop.\ \propdeltaperbel). At iteration $k$, we have a $p$-stable policy $\m^k$, and we generate a $p$-stable policy $\m^{k+1}$ according to
$$\m^{k+1}(x)\in \argmin_{u\in U(x)}\big\{g(x,u)+J_{\m^k,\d_k}\big(f(x,u)\big)\big\},\qquad x\in X .\xdef\intermrelz{\lab}\eqnum\show{oneo}$$
Note that under Assumption \assumptionpertpi, the algorithm is well-defined, and 
is guaranteed to generate a sequence of $p$-stable stationary policies. 

We will use for all policies $\m$ and scalars $\d>0$ the mappings  $T_\m:{\cal E}^+(X)\mapsto {\cal E}^+(X)$ and $T_{\m,\d}:{\cal E}^+(X)\mapsto {\cal E}^+(X)$ by
$$(T_\m J)(x)=g\big(x,\m(x)\big)+J\big(f(x,\m(x))\big),\qquad x\in X,$$
$$(T_{\m,\d} J)(x)=g\big(x,\m(x)\big)+\d p(x)+J\big(f(x,\m(x))\big),\qquad x\in X,$$
and the mapping $T:{\cal E}^+(X)\mapsto {\cal E}^+(X)$ given by
$$(TJ)(x)=\inf_{u\in U(x)}\big\{g(x,u)+J\big(f(x,u)\big)\big\},\qquad x\in X.$$
For any integer $m\ge1$, we denote by $T_\m^m$ and $T_{\m,\d}^m$ the $m$-fold compositions of the mappings $T_\m$ and $T_{\m,\d}$, respectively. We have the following proposition.

\xdef\proppertpi{\propn}\propnum\show{myproposition}

\texshopbox{\proposition{\proppertpi:} Let Assumption \assumptionpertpi\ hold. Then  for a sequence of $p$-stable policies $\{\m^k\}$ generated by the perturbed PI algorithm \intermrelz, we have $J_{\m^k,\d_k}\downarrow\hat J_p$ and $J_{\m^k}\to\hat J_p$.}

\proof The algorithm definition \intermrelz\ implies that for all integer $m\ge1$ we have for all  $x_0\in X$,
$$J_{\m^k,\d_k}(x_0)\ge (TJ_{\m^k,\d_k})(x_0)+\d_k p(x_0)=   (T_{\m^{k+1},\d_k}J_{\m^k,\d_k})(x_0)\ge (T_{\m^{k+1},\d_k}^mJ_{\m^k,\d_k})(x_0)\ge (T_{\m^{k+1},\d_k}^m\bar J)(x_0),$$
where $\bar J$ is the identically zero function [$\bar J(x)\equiv0$]. From this relation we obtain 
$$J_{\m^k,\d_k}(x_0)\ge \lim_{m\to\infty}(T_{\m^{k+1},\d_k}^m\bar J)(x_0)=\lim_{m\to\infty}\lf\{\sum_{\ell=0}^{m-1}\big(g(x_\ell,\m^{k+1}(x_\ell))+\d_k p(x_\ell)\big)\ri\}\ge J_{\m^{k+1},\d_{k+1}}(x_0),$$
as well as
$$J_{\m^k,\d_k}(x_0)\ge (TJ_{\m^k,\d_k})(x_0)+\d_k p(x_0)\ge J_{\m^{k+1},\d_{k+1}}(x_0).$$
It follows that $\{J_{\m^k,\d_k}\}$ is monotonically nonincreasing, so that $J_{\m^k,\d_k}\downarrow J_\infty$ for some $J_\infty$, and 
$$\lim_{k\to\infty}TJ_{\m^k,\d_k}= J_\infty.\xdef\intermlimit{\lab}\eqnum\show{oneo}$$ 

 We also have, using the fact $J_\infty\le J_{\m^k,\d_k}$, 
$$\eqalignno{\inf_{u\in U(x)}\big\{g(x,u)+J_\infty\big(f(x,u)\big)\big\}&\le\lim_{k\to\infty}\inf_{u\in U(x)}\big\{g(x,u)+J_{\m^k,\d_k}\big(f(x,u)\big)\big\}\cr
&\le\inf_{u\in U(x)}\lim_{k\to\infty}\big\{g(x,u)+J_{\m^k,\d_k}\big(f(x,u)\big)\big\}\cr
&=\inf_{u\in U(x)}\lf\{g(x,u)+\lim_{k\to\infty}J_{\m^k,\d_k}\big(f(x,u)\big)\ri\}\cr
&= \inf_{u\in U(x)}\big\{g(x,u)+J_\infty\big(f(x,u)\big)\big\}.\cr}$$
Thus equality holds throughout above, so that
$$\lim_{k\to\infty}TJ_{\m^k,\d_k}= TJ_\infty.$$
Combining this with Eq.\ \intermlimit, we obtain $ J_\infty= TJ_\infty$, i.e., $J_\infty$ solves Bellman's equation. We also note that $J_\infty\le J_{\m^0,\d_0}$ and that $J_{\m^0,\d_0}\in S_p$ by assumption, so that $J_\infty\in S_p$.
By Prop.\ \proptaffmonft(a), it follows that $J_\infty=\hat J_p$. \qed

Note that despite the fact $J_{\m^k}\to\hat J_p$, the generated sequence $\{\m^k\}$ may exhibit some serious pathologies in the limit. In particular, if $U$ is a metric space and $\{\m^k\}_{\cal K}$ is a subsequence of policies that converges to some $\bar \m$, in the sense that 
$$\lim_{k\to\infty,\,k\in{\cal K}}\m^k(x)=\bar \m(x),\qquad\forall\ x\in X,$$
it does not follow that $\bar \m$ is $p$-stable. In fact it is possible to construct examples where the generated sequence of $p$-stable policies $\{\m^k\}$ satisfies $\lim_{k\to\infty}J_{\m^k}=\hat J_p=\jstar$, yet $\{\m^k\}$ may converge to a $p$-unstable policy whose cost function is strictly larger than $\hat J_p$. Example 2.1 of the paper [BeY16] provides an instance of a stochastic  shortest path problem with two states, in addition to the termination state, where this occurs.

\subsection{An Optimistic Policy Iteration Method}

\pn Let us  consider an optimistic variant of PI, where policies are evaluated inexactly, with a finite number of VIs. We use a fixed forcing function $p$. We will show that the algorithm can be used to compute $\hat J_p$, the restricted optimal cost function over the $p$-stable policies. The algorithm generates a sequence $\{J_k,\m^k\}$ according to
$$T_{\m^k} J_k=T J_k,\qquad J_{k+1}=T_{\m^k}^{m_k}J_k,\qquad k=0,1,\ldots,\xdef\optpoiter{\lab}\eqnum\show{oneo}$$
where $m_k$ is a positive integer for each $k$. We assume that a policy $\m^k$ satisfying $T_{\m^k} J_k=T J_k$ can be found for all $k$, but it need not be $p$-stable. However, the algorithm requires that 
$$J_0\in{\cal J},\qquad J_0\ge TJ_0,\qquad J_0\in{\cal W}_p.\xdef\optpicond{\lab}\eqnum\show{oneo}$$
 This may be a restrictive assumption. 
We have the following proposition. 

\xdef\propdetpolicyitopt{\propn}\propnum\show{myproposition}
\texshopbox{\proposition{\propdetpolicyitopt: (Convergence of Optimistic PI)}Assume that there exists at least one $p$-stable policy $\p\in\P_p$, and that $J_0$ satisfies Eq.\ \optpicond. Then a sequence $\{J_k\}$ generated by the optimistic PI algorithm \optpoiter\ belongs to ${\cal W}_p$ and satisfies $J_{k}\downarrow \hat J_p$.
}
 \proof Since $J_0\ge \hat J_p$ and $\hat J_p=T\hat J_p$ [cf.\ Prop.\  \propbeleqterm(a)], all operations on any of the functions $J_k$ with $T_{\m^k}$ or $T$ maintain the inequality $J_k\ge \hat J_p$ for all $k$, so that $J_k\in{\cal W}_p$ for all $k$. Also the conditions $J_0\ge TJ_0$ and  $T_{\m^k} J_k=T J_k$ imply that
 $$J_0= J_1\ge T_{\m^0}^{m_0+1}J_0=T_{\m^0}J_1\ge TJ_1=T_{\m^1}J_1\ge\cdots\ge J_2,\eqnum\show{oneo}$$
and continuing similarly, 
$$J_k\ge TJ_k\ge J_{k+1},\qquad k=0,1,\ldots.\xdef\bracketineq{\lab}\eqnum\show{oneo}$$
Thus  $J_{k}\downarrow J_\infty$ for some $J_\infty$, which must satisfy $J_\infty\ge \hat J_p$, and hence belong to $ {\cal W}_p$. By taking limit as $k\to\infty$ in Eq.\ \bracketineq\ and using an argument similar to the one in the proof of Prop.\ \proppertpi, it follows that $J_\infty=TJ_\infty$. By Prop.\  \propbeleqterm(a), this implies that $J_\infty\le \hat J_p$. Together with the inequality $J_\infty\ge \hat J_p$ shown earlier, this proves that $J_\infty=\hat J_p$.
\qed

As an example,  for the shortest path problem of Example \examplecounterpi, the reader may verify that for the case where $p(x)=1$, for $x=1$, the optimistic PI algorithm converges in a single iteration to 
$$\hat J_p(x)=\cases{1&if $x=1$,\cr
0&if $x=0$,\cr}$$
provided that
$J_0\in  {\cal W}_p=\big\{J\mid J(1)\ge 1,\ J(0)=0\big\}.$
For other starting functions $J_0$, the algorithm converges in a single iteration to the function
$$J_\infty(1)=\min\big\{1,\,J_0(1)\big\},\qquad J_\infty(0)=0.$$
All functions $J_\infty$ of the form above are solutions of Bellman's equation, but only $\hat J_p$ is restricted optimal.

\vskip-1.0pc

\section{Concluding Remarks}

\pn We have considered deterministic optimal control problems with a cost-free and absorbing destination under general assumptions, which include arbitrary state and control spaces, and a Bellman's equation with multiple solutions. Within this context, we have used perturbations of the cost per stage and the ideas of semicontractive DP as a means to connect classical issues of stability and optimization. In particular, we have shown that the restricted optimal cost function over just the stable policies is a solution of Bellman's equation, and that versions of the VI and PI algorithm are attracted to it. Moreover, the restricted optimal cost $J^+$ over the ``fastest" policies (the ones that terminate) is the largest solution of Bellman's equation. The generality of our framework makes our results a convenient starting point for analysis of related problems and methods, involving additional assumptions,  and/or cost function approximation and state space discretization. 

An interesting open question is how to discretize continuous-spaces problems to solve Bellman's equation numerically. As an example, consider the linear-quadratic problem of Example \examplelqminenergy. Any reasonable discretization of this problem is a finite-state (deterministic or stochastic) shortest path problem, whose Bellman equation has a unique solution that approximates the solution $\jstar$ of the continuous-spaces problem, while missing entirely the solution $J^+$. The same is true for the optimal stopping problem of Example \exampleperturb. In such cases, one may discretize a $\d$-perturbed version of the problem, which is better behaved, and use a small value of $\d$ to obtain an approximation to $J^+$. However, the limiting issues as $\d\downarrow0$ remain to be explored.

\vskip-1pc

\section{References}
\vskip-0.9pc
\def\ref{\vskip1.pt\pn}

\ref[AnM07] Anderson, B.\ D., and Moore, J.\ B., 2007.\ Optimal Control: Linear Quadratic Methods, Courier Corporation.

\ref [BeS78] Bertsekas, D.\ P., and Shreve, S.\ E., 1978.\ Stochastic Optimal Control: The Discrete Time Case. New York: Academic Press; may be downloaded from http://web.mit.edu/dimitrib/www/home.html.

\ref [BeT96]  Bertsekas, D.\ P., and Tsitsiklis, J.\ N., 1996.\ Neuro-Dynamic
Programming, Athena Scientific, Belmont, MA.

\ref[BeY16] Bertsekas, D.\ P., and Yu, H., 2016.\ ``Stochastic Shortest Path Problems Under Weak Conditions,"  Lab.\ for Information and Decision Systems Report LIDS-2909. 

\ref[Ber12] Bertsekas, D.\ P., 2012.\ Dynamic Programming and Optimal Control, Vol.\ II: Approximate Dynamic Programming, Athena Scientific, Belmont, MA.

\ref[Ber13] Bertsekas, D.\ P., 2013.\ Abstract Dynamic Programming, Athena Scientific, Belmont, MA; a second edition appeared in 2017 on-line at http://web.mit.edu/dimitrib/www/home.html.

\ref[Ber14] Bertsekas, D.\ P., 2014.\ ``Robust Shortest Path Planning and Semicontractive Dynamic Programming," Lab.\ for Information and Decision Systems Report LIDS-P-2915, MIT, Feb. 2014 (revised Jan. 2015 and June 2016); arXiv preprint arXiv:1608.01670; to appear in Naval Research Logistics.

\ref[Ber15a] Bertsekas, D.\ P., 2015.\ ``Value and Policy Iteration in Deterministic Optimal Control and Adaptive Dynamic Programming," arXiv preprint arXiv:1507.01026; IEEE Trans.\  on Neural Networks and Learning Systems, Vol.\ 28, 2017, pp.\ 500-509.

\ref[Ber15b] Bertsekas, D.\ P., 2015.\ ``Regular Policies in Abstract Dynamic Programming," Lab. for Information and Decision Systems Report LIDS-P-3173, MIT; arXiv preprint arXiv:1609.03115; to appear in SIAM J.\ on Optimization.

\ref[Ber16] Bertsekas, D.\ P., 2016.\ ``Affine Monotonic and Risk-Sensitive Models in Dynamic Programming", Lab. for Information and Decision Systems Report LIDS-3204, MIT, June 2016; arXiv preprint arXiv:1608.01393.

\ref[Ber17] Bertsekas, D.\ P., 2017.\ Dynamic Programming and Optimal Control, Vol.\ I, 4th edition, Athena Scientific, Belmont, MA.

\ref[Hey14] Heydari, A., 2014.\ ``Stabilizing Value Iteration With and Without Approximation Errors," available at arXiv:1412.5675.

\ref[JiJ14] Jiang, Y., and Jiang, Z.\ P., 2014.\ ``Robust Adaptive Dynamic Programming and Feedback Stabilization of Nonlinear Systems," IEEE Trans.\ on Neural Networks and Learning Systems, Vol.\ 25, pp.\ 882-893.

\ref[LLL08] Lewis, F.\ L., Liu, D., and Lendaris, G.\ G., 2008.\ Special Issue on Adaptive Dynamic Programming and Reinforcement Learning in Feedback Control, IEEE Trans.\ on Systems, Man, and Cybernetics, Part B, Vol.\ 38, No.\ 4.

\ref[LeL13] Lewis, F.\ L., and Liu, D., (Eds), 2013.\ Reinforcement Learning and Approximate Dynamic Programming for Feedback Control, Wiley, Hoboken, N.\ J.

\ref[LiW13] Liu, D., and Wei, Q., 2013.\ ``Finite-Approximation-Error-Based Optimal Control Approach for Discrete-Time Nonlinear Systems, IEEE Trans.\  on Cybernetics, Vol.\ 43, pp.\ 779-789.

\ref [Kle68] Kleinman, D.\ L., 1968.\  ``On an Iterative Technique for Riccati
Equation Computations," IEEE Trans.\ Aut.\  Control, Vol.\ AC-13, pp.\ 114-115.

\ref[Kuc72] Kucera, V., 1972.\ ``The Discrete Riccati Equation of Optimal Control," Kybernetika, Vol.\ 8, pp.\ 430-447.

\ref[Kuc73] Kucera, V., 1973.\ ``A Review of the Matrix Riccati Equation," Kybernetika, Vol.\ 9, pp.\ 42-61.

\ref[LaR95] Lancaster, P., and Rodman, L., 1995.\ Algebraic Riccati Equations, Clarendon Press, Oxford, UK.

\ref[Put94] Puterman, M.\ L., 1994.\ Markov Decision Processes: Discrete Stochastic Dynamic Programming, Wiley, N.\ Y. 

\ref[SBP04] Si, J., Barto, A., Powell, W., and Wunsch, W., (Eds.), 2004.\ Learning and Approximate Dynamic
Programming, IEEE Press, N.\ Y.

\ref [Str66] Strauch, R., 1966.\  ``Negative Dynamic Programming," Ann.\ Math.\
Statist., Vol.\ 37, pp.\ 871-890.

\ref[VVL13] Vrabie, D., Vamvoudakis, K.\ G., and Lewis, F.\ L., 2013.\ Optimal Adaptive Control and Differential Games by Reinforcement Learning Principles,
The Institution of Engineering and Technology, London.

\ref[Wil71] Willems, J., 1971.\ ``Least Squares Stationary Optimal Control and the Algebraic Riccati Equation," IEEE Trans.\ on Automatic Control, Vol.\ 16, pp.\ 621-634.

\end